\documentclass[twoside, 11pt]{article}
\usepackage{mathrsfs,amsfonts,amsmath}
\RequirePackage{ifthen,calc}
\usepackage{theorem}
\usepackage[dvips]{color}
\usepackage{graphicx}

\def\R{{\mathbb R}}

\def\Re {{\rm Re}\,}
\def\E{{\mathbb E}}
\def\P{{\mathbb P}}

\def\C{{\mathbb C}}
\def\L{{L^{p}(T)}}

 \catcode`@=11
 \def\@evenhead{\hbox to\textwidth{\footnotesize\rm\thepage \hfill
  {\it }}} 

 \def\@oddhead{\hbox to \textwidth{\footnotesize{\it
 } \hfill\thepage}}

\renewcommand{\subsection}{\makeatletter
 \renewcommand{\@seccntformat}[1]{{\csname the##1\endcsname.}\hspace{0.45em}}
 \makeatother \@startsection
{subsection}
{1}
{0pt}
{\baselineskip}
{0.5\baselineskip}
{\normalsize\bfseries\mathversion{bold}}}

\catcode`@=12

\newtheorem{thm}{\noindent Theorem}[section]
\newtheorem{lem}{\noindent Lemma}[section]
\newtheorem{cor}{\noindent Corollary}[section]

\theoremheaderfont{\normalfont\bfseries}
\theorembodyfont{\slshape} {\theorembodyfont{\rmfamily}} {\theorembodyfont{\rmfamily}
\newtheorem{defn}{\noindent Definition}[section]}
{\theorembodyfont{\rmfamily} } {\theorembodyfont{\rmfamily}
}
{\theorembodyfont{\rmfamily} } {\theorembodyfont{\rmfamily}
\newtheorem{rem}{\noindent Remark}[section]}
{\theorembodyfont{\rmfamily} } {\theorembodyfont{\rmfamily} } {\theorembodyfont{\rmfamily}
}

 \def\beqlb{\begin{eqnarray}}\def\eeqlb{\end{eqnarray}}
 \def\beqnn{\begin{eqnarray*}}\def\eeqnn{\end{eqnarray*}}
 \def\Re{\textrm{\it Re}}
 \numberwithin{equation}{section}
\def\qed{\hfill$\square$\smallskip}

\begin{document}
\title{Exact tail asymptotics for a three dimensional Brownian-driven tandem queue with  intermediate inputs\footnotetext{\hspace{-5ex}
${[1]}$ School of Statistics, Shandong University of Finance and Economics,
Jinan 250014,  China.
\\${[2]}$ School of Mathematics, Carleton University, Ottawa ON  K1S5B6, CA
\newline
}}
\author{\small Hongshuai Dai$^{[1,2]}$, Dawson A. Donald$^{[2]}$, Yiqiang Q. Zhao$^{[2]}$ }
\maketitle
\begin{abstract}

The semimartingale reflecting Brownian motion (SRBM) can be a heavy traffic limit for many server queueing networks. Asymptotic properties for stationary probabilities of the SRBM have attracted a lot of attention recently.  However, many results are obtained only for the two-diemnsional SRBM.  There is only little work related to higher dimensional ($\geq 3$) SRBMs.  In this paper, we  consider a three dimensional SRBM: A three dimensional Brownian-driven tandem queue with  intermediate inputs. We are interested in tail asymptotics for stationary distributions.  By generalizing the kernel method and using coupla, we obtain  exact tail asymptotics for the marginal stationary distribution of  the buffer content in the third buffer and the joint stationary distribution.
\end{abstract}

\small {{\bf MSC(2000):}  60K25, 60J10.

}
\section{Introduction}
Since  Harrison and Reiman \cite{HR1981a,HR1981b}, Varadhan and Williams\cite{VW1985}, and Williams\cite{W1985a,W1985b} introduced the semimartingale reflecting Brownian motion, SRBM has received a lot of attention. Stationary properties of stationary distributions of SRBM when they exist, are  important, especially in applications. However, except for a very limited number of special cases, a simple closed expression for the stationary distribution is not available. Hence, exact tail behaviour of stationary distributions becomes most important. Recently, many results about two-dimensional SRBM have been obtained. Harrison and Hasenbein~\cite{HH2009} presented sufficient and necessary conditions for the existence of a stationary distribution.  Dai and Miyazawa \cite{DM2011} studied exact tail asymptotics for the marginal distributions of SRBM by using a geometric method.  Dai, Dawson and Zhao \cite{DDZ} applied the kernel method to obtain exact tail asymptotics for the boundary measures of SRBM.  Franceschi and Kurkova \cite{FK2016} studied exact tail asymptotics of the stationary distribution along some path by analytic methods. Franceschi and Raschel \cite{FR2017} studied exact tail behaviour of the boundary stationary distributions of SRBM  by using the boundary value problems. However,  we note that all aforementioned results are only for two-dimensional SRBM. In this paper, we will  consider a three dimensional SRBM.

Miyazawa and Rolski \cite{MR2009} generalized the result of Lieshout and Mandjes \cite{LM2007, LM2008}, and studied a two-dimensional L\'evy-driven tandem queue with an intermediate input.  They obtained exact tail asymptotics for the Brownian inputs, while weaker tail asymptotic results  were obtained for the general L\'evy input. They also tried to discuss  higher-dimensional cases. However  only  the stationary equation was obtained in terms of moment generating functions, and  tail asymptotic properties for the marginal distributions were left for a future work.  In this paper, we consider a three-dimensional Brownian-driven tandem queue with intermediate inputs.  We  derive exact tail aysmptotics for the marginal stationary distribution of the third buffer content,  since exact tail asymptotic results for the first two buffer content can be obtained directly from results for two-dimensional SRBM.  We also note that all results related to exact tail asymptotics for stationary distributions of SRBM are only for marginal stationary distributions and boundary stationary distributions. There are no results referred in the literature on asymptotic properties for the joint stationary distribution of SRBM, which is also considered in this paper.

 In this paper, we apply both the Kernel method and extreme value theory to study  tail asymptotics. The kernel method  has been systematically applied to study random walks in the quarter plane by Li and Zhao \cite{LZ2012} and references therein.
Key steps in applying the kernel method for random walks in the quarter plane are:(i) Establishing the fundamental form:
\[
    h(x,y) \pi(x,y) = h_1(x,y) \pi_1(x) + h_2(x,y) \pi_2(y) +
    h_0(x,y) \pi_{0,0},
\]
where $\pi(x,y)$, $\pi_1(x)$ and $\pi_2(y)$ are unknown generating
functions for joint and two boundary probabilities, respectively.
(ii)  Finding a branch
$Y=Y_0(x)$ such that $h(x,Y_0(x))=0$, which leads to
a relationship between the two unknown boundary generating
functions:
\beqlb\label{q-1}
    h_1(x,Y_0(x)) \pi_1(x) + h_2(x,Y_0(x)) \pi_2(Y_0(x)) +
    h_0(x,Y_0(x)) \pi_{0,0}=0;
\eeqlb
(iii) Based on \eqref{q-1},  carrying out a singularity analysis for $\pi_1$ and
$\pi_2$, which leads to not only a decay rate, but also exact tail asymptotic properties of the boundary probabilities through a Tauberian-like theorem.  In this paper, we will extend this method to study a three-dimensional SRBM.   By using the kernel method, we can get exact tail asymptotics for the marginal stationary distributions.

In this paper, we also study asymptotic properties for the joint stationary distributions. However, we cannot use the kernel method to study tail behaviours of the joint stationary distribution, since the kernel method relies on the Tauberian-like Theorem, which is valid only for univariate functions.  By using the kernel method, we can  get tail equivalence for the marginal distributions, from which we will further study the tail dependence of the joint stationary distribution. Tail dependence describes the amount of dependence in the upper tail or lower tail of a multivariate distributions and has been widely used in extreme value analysis and in quantitative risk management. Once we get the dependence,  we can study multivariate extreme value distribution of the joint stationary distribution. The extreme value distribution is very useful since from a sample of vectors of maximum, one can make inferences about the upper tail of the stationary distribution using multivariate extreme value theory. Based on the multivariate extreme distribution, by using copula, we can get tail behaviour of the joint stationary distributions.

In this paper, we  study a three dimensional SRBM and anticipate the tools developed in this paper will  be useful in analyzing the general $d$-dimensional case. The rest of this paper is organized as follows:  In Section 2, a three-dimensional Brownian-driven tandem queue with  intermediate inputs is introduced.  To apply the kernel method for asymptotic properties for the marginal $L_3$, we study the kernel equation and the analytic continuation of moment generating functions in Section 3. We study some asymptotic properties of  moment generating functions in Section 4.  Asymptotic results for the marginal distributions are present in Section 5. In Section 6, we study asymptotic properties of the joint stationary distribution.

\section{Model and Preliminaries}

In this section, we  introduce a three dimensional Brownian-driven tandem queue with intermediate inputs and establish a stationary equation satisfied by stationary probabilities.  This tandem queue has three nodes, numbered as 1,2, 3, each of which has exogenous input process and a constant processing rate.   Outflow from the node 1 goes to  node 2,  and the outflow from node 2 goes to  node 3. Finally, outflow from  node $3$ leaves the system,  see Fig.1 below.

\begin{center}\includegraphics[width=9cm,height=1.5cm]{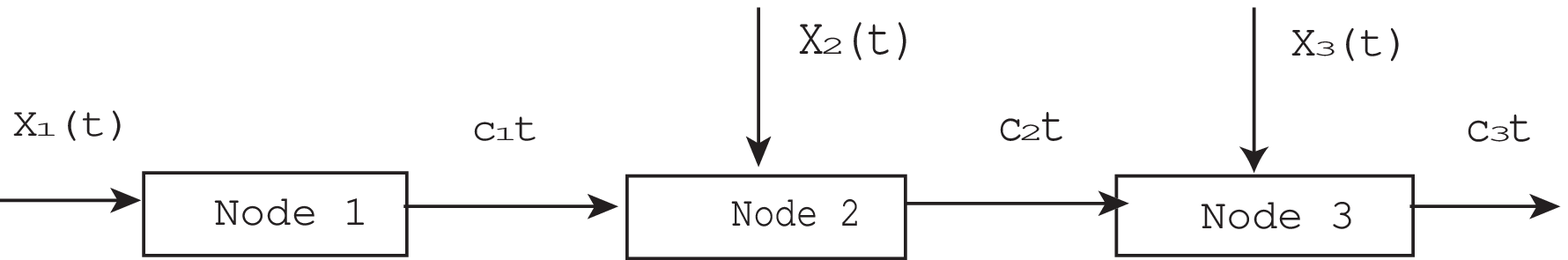}\\
\vspace{3mm}
{\footnotesize {\bf Fig. 1} A  tandem queue with 3 nodes.}
\end{center}

 We assume that the exogenous inputs are Brownian processes of the form:
 \beqlb\label{s2-1}
 X_i(t)=\lambda_i t+B_i(t), \; i=1,2,3,
 \eeqlb
where $\lambda_i>0$  is a nonnegative constant, and $B_i(t)$ is a Brownian motion with variance $\sigma^2_i$ and no drift. Without loss of generality, we assume that the correlation coefficients $\rho_{B_iB_j}<1$, $i,j=1,2,3$, where $B(1)=(B_1,B_2,B_3)'$. Denote the processing rate at  node $i$ by $c_i>0$. Let $L_i(t)$ be the buffer content at node $i$ at time $t\geq 0$ for $i=1,2,3$, which are formally defined as
\beqlb\label{s2-2}
L_1(t)&&=L_1(0)+X_1(t)-c_1t+Y_1(t),
\label{s2-3}
\\L_i(t)&&=L_i(0)+X_i(t)+c_{i-1}t-c_{i} t-Y_{i-1}(t)+Y_i(t),\; i\geq 2,
\eeqlb
where $Y_i(t)$ is a regulator at  node $i$, that is,  a minimal nondecreasing process for $L_i(t)$ to be nonnegative.  In fact, we can  regard $\big(L_1(t),L_2(t),L_3(t)\big)'$ as a reflection mapping from the net flow processes $\big(X_1(t)-c_1t, X_2(t)-c_2 t, X_3(t)+c_{2} t-c_3 t\big)'$ with the reflection matrix
\beqlb\label{s2-5}
R=\left[ {\begin{array}{*{20}{c}}
{{1}}\quad & {{0}}\quad& {{0}} \\
{{-1}} \quad& {{1}}\quad & {{0}} \\
{{0}} \quad& {{-1}}\quad & {{1}} \\
\end{array}} \right]=[R^1,R^2,R^3].
\eeqlb
Let  $L(t)=\big(L_1(t),L_2(t),L_3(t)\big)'$ and $B(t)=\big(B_1(t),B_2(t),B_3(t)\big)'$. Then
\beqnn
L(t)=B(t)+\Lambda t+R Y(t)+L(0),
\eeqnn
where $\Lambda=\big(\lambda_1-c_1,\lambda_2+c_{1}-c_2,\lambda_3+c_{2}-c_3\big)'$ and $Y(t)=\big(Y_1(t),Y_2(t),Y_3(t)\big)'$.

Without any difficulty,  we can obtain that the tandem queue has the stationary distribution if and only if
\beqlb\label{s2-6a} \sum_{i=1}^k \lambda_i\leq c_k,\;1\leq k\leq 3.\eeqlb
Moreover, by Harrison and Williams \cite{HW1987}, we can  get that the stationary distribution of $\{L(t)\}$ is unique. Throughout this paper, we denote this stationary distribution by $\pi$. In order to simplify the discussion, in this paper, we  refine the stability condition \eqref{s2-6a} to assume that
\beqlb\label{s2-6}
\lambda_1<c_1,\; \textrm{and}\;\lambda_{i}+c_{i-1}<c_i,\; i=2,3.
\eeqlb
\begin{rem}
From the  proofs of the main results of this paper, it is clear that under the more general stability condition \eqref{s2-6a}, we can use the same argument to discuss  tail asymptotics. The only difference  is that we need to discuss possible relationships between the parameters $\lambda_i$ and $c_i$, $i=1,2,3$, before  we use the arguments in the proofs in this paper.  For each of the possible relationships,  we repeat  the method applied in this paper to study  tail asymptotics.
\end{rem}

We are interested in  asymptotic tail behaviour of the stationary distribution.  Recall that  a positive function $g(x)$ is said to have exact tail asymptotic $h(x)$, if
\beqnn
\lim_{x\to\infty}\frac{g(x)}{h(x)}=1.
\eeqnn
In this paper, our main aim is to find  exact tail asymptotics for various stationary distributions. Moment generating function will play an important role in determining these exact tail asymptotics.  We first introduce  moment generating functions for stationary distributions.  Let $L=(L_1,L_2,L_3)'$ be the stationary random vector with the stationary distribution $\pi$.  The moment generating function $\phi(\cdot)$ for $L$ is given  by:
\beqlb \label{s2-9}
\phi(x,y,z)=\E \big[e^{ x L_1+ y L_2+ z L_3}\big],\;\textrm{for any}\;(x,y,z)'\in \mathbb{R}^3.
\eeqlb

We apply the kernel method to study  tail asymptotics for stationary distributions.  In order to apply the kernel method, we need establish a relationship between the moment generating function $\phi(\cdot)$ for the stationary distribution and the moment generating functions for the boundary measures defined below. For any Borel set $A\subset \mathscr{B}(\R^3)$, we define the boundary measures $V_i(\cdot)$,  $i=1,2,3$, by
\beqlb\label{s2-8}
V_i(A)=\E_\pi\Big[\int_0^1 I_{\{L(u)\in A\}}d Y_i(u)\Big].
\eeqlb
Moreover, due to Harrison and Williams \cite{HW1987}, we obtain that the density functions for $V_i$, $i=1,2,3,$ exist. Then, their  moment generating functions are defined by
\beqlb \label{s2-10}\phi_i(x,y,z)=\int_{\R_+^3}e^{<\vec{x},\theta>}V_i(d\theta)=\E\Big[\int_0^1 e^{<\vec{x},\;L(t)>}dY_i\Big], \;i=1,2,3,\eeqlb
where $\vec{x}=(x,y,z)'\in\mathbb{R}^3$.

Next, we establish the relationship between these moment generating functions. In fact,
there is a nice connection during them. The following lemma is due to Konstantopoulos, Last and Lin \cite{KLL2004}.
\begin{lem}\cite[Theorem 4]{KLL2004}\label{lem-a-7}  For each $(x,y,z)'\in\R^3$ with $\phi(x,y,z)<\infty$ and $\phi_i(x,y,z)<\infty$, $i=1,2,3$, we have
\beqlb\label{s2-69}
H(x,y,z)\phi(x,y,z)=H_1(x,y)\phi_1(x,y,z)+H_2(y,z)\phi_2(x,y,z)+H_3(z)\phi_3(x,y,z),
\eeqlb
where
\beqlb\label{s2-70}
&&H(x,y,z)=-\frac{1}{2}\big(\sigma_1^2x^2+\sigma_2^2y^2
+\sigma_3^2z^2\big)+(c_1-\lambda_1)x+(c_2-\lambda_2-c_{1})y+(c_3-\lambda_3-c_{2})z,\qquad
\label{s2-71}
\\&&H_1(x,y)=x-y,
\label{s2-72}
\\&&H_2(y,z)=y-z,
\label{s2-73}
\\&&H_3(z)=z.
\eeqlb
\end{lem}

From Lemma \ref{lem-a-7}, we can prove the following lemma.
\begin{lem}\label{s2-lem2} For $1\leq k\leq 3$, we have
\beqlb
\label{s2-85}
\phi_k(0,0,0)&&=\E_{\pi}[Y_k(1)]=c_k-\sum_{i=1}^k\lambda_i>0.
\eeqlb
\end{lem}
{\it Proof:}  From \eqref{s2-10}, we get that
\eqref{s2-69} makes sense for $x\in\R^3\setminus\R_{+}^3$.  Let $x=(0,x_k,0)'$ with $x_k< 0$. From \eqref{s2-69}, we get
\beqnn
&&\big((c_k-\lambda_k-c_{k-1})x_k-\frac{1}{2}\sigma_k^2 x_k^2\big)\phi(0,x_k,0)
=-x_{k}\phi_{k-1}(0,x_k,0)+x_k\phi_k(0,x_k,0),
\eeqnn
i.e.,
\beqlb\label{s2-30}
&&\big((c_k-\lambda_k-c_{k-1})-\frac{1}{2}\sigma_k^2 x_k\big)\phi(0,x_k,0)\nonumber
=-\phi_{k-1}(0,x_k,0)+\phi_k(0,x_k,0),
\eeqlb
where $c_0=0$ and $\phi_0=0$.
Letting $x_k$ go to $0$ in \eqref{s2-30}, we get that the left-hand side of  equation \eqref{s2-30} equals to
\beqlb\label{s2-30-1}
c_k-\lambda_k-c_{k-1}
\eeqlb
since $\phi(0,0,0)=1$.  Hence,
\beqlb\label{s2-30-2}
(c_k-\lambda_k-c_{k-1})+\phi_{k-1}(0,0,0)=\phi_k(0,0,0).
\eeqlb
Let $k=1$. Then, one can easily get that
\beqlb\label{s2-30-3}
\phi_1(0,0,0)=c_1-\lambda_1.
\eeqlb
By \eqref{s2-30-2} and \eqref{s2-30-3}, we can get the lemma holds.
 \qed

In general, it is difficult or impossible to obtain the explicit expression for the stationary distribution $\pi$ or its moment generating function.  However in some special cases, it becomes possible. For example, if there are no intermediate inputs, that is $X_k=0$, $k=2,3$, Miyazawa and Rolski \cite{MR2009} obtained an explicit expression of $\phi(\theta)$. For a general case, our focus is on its tail asymptotics.  There are a few aviable methods for studying tail asymptotics, for example, in terms of large deviations and boundary value problems.  In this paper, we study  tail asymptotics of the marginal distribution $\P\big(L_3<x\big)$ via the kernel method introduced by Li and Zhao \cite{LZ2011} and  asymptotic properties of the joint stationary distribution by extreme value theory and copula.

At the end of this section, we present a technical lemma, which  plays an important role in finding the tail asymptotics of the marginal distribution $\mathbb{P}(L_3<x)$.
\begin{lem}\label{s2-lem1}
$\phi(0,0,x)$ and $\phi_{2}(0,0,x)$ have the same singularities.
\end{lem}
{\it Proof:}
Let $\theta=(0,0,x)'$. Then,
\beqlb\label{s2-17}
H(0,0,x)=-\frac{1}{2}\sigma_3^2 x^2+(c_3-\lambda_3-c_{2})x.
\eeqlb
Note that for any $(x,y,z)'$, since $\Delta Y_3(t)>0$ only if $L_3(t)=0$, we have
\beqnn
\phi_3(x,y,z)=\phi_3(x,y,0).
\eeqnn
Then, by (\ref{s2-69}) and (\ref{s2-85}),
\beqlb\label{s2-18}
H(0,0,x)\phi(0,0,x)=-x\phi_{2}(0,0,x)+x\big(c_3-\sum_{i=1}^3\lambda_i\big).
\eeqlb
From \eqref{s4-97} below, we  get that $$\phi_{2}(0,0,\frac{2(c_3-\lambda_3-c_{2})}{\sigma_3^2})<+\infty.$$  Letting
$x=\frac{2(c_3-\lambda_3-c_{2})}{\sigma_3^2}$ in (\ref{s2-18}), we obtain
\beqnn\label{s2-19}
-\phi_{2}(0,0,\frac{2(c_3-\lambda_3-c_{2})}{\sigma_3^2})+\big(c_3-\sum_{i=1}^3\lambda_i\big)=0,
\eeqnn
i.e.,
\beqlb\label{s2-20}
\phi_{2}(0, \frac{2(c_3-\lambda_3-c_{2})}{\sigma_3^2})=c_3-\sum_{i=1}^3 \lambda_i.
\eeqlb
Therefore, by (\ref{s2-18}) and (\ref{s2-20}),
\beqlb\label{s2-21}
\phi(0,0,x)=\frac{\phi_{2}\big(0,0,\frac{2(c_3-\lambda_3-c_{2})}{\sigma_3^2}\big)-\phi_{2}(0,0,x)}{-\frac{1}{2}\sigma_3^2x+(c_3-\lambda_3-c_{2})}.
\eeqlb
By (\ref{s2-20}) and (\ref{s2-21}), one can easily get that $x=\frac{2(c_3-\lambda_3-c_{2})}{\sigma_3^2}$ is a removable singularity of $\phi(0,0,x)$. The proof of this lemma is completed.\qed

\section{Kernel Equation and Analytic Continuation}
\label{S3}
In this paper, we  apply the kernel method to study  tail asymptotics for the marginal stationary measure $\P(L_3<x)$.  In order to do it,  we need  the Tauberian-like Theorem (Theorem \ref{TLT}). For applying this theorem,  we need to study the analytic properties of the moment generating function $\phi_2(x,y,z)$.
\subsection{Kernel Equation and Branch Points}
To study analyitic properties of the moment generating functions, we first  focus on the kernel equation and the corresponding branch points. For this purpose, we consider the kernel equation:
\beqlb\label{s3-1}
H(x,y,z)=0,
\eeqlb
which is critical in our analysis.

Since tail asymptotics for $\mathbb{P}(L_3<z)$ is our focus, we first treat $z$ in $ (x,y,z)'\in \R^3$ as a variable.  Inspired by the procedure of applying the kernel method, for example, see Li and Zhao \cite{LZ2011},  we first construct the  relationship between $z$ and  $x$, $y$. The kernel equation in \eqref{s3-1} defines an implicit function  $z$ in variables $x$ and $y$ when we only consider non-negative values for $z$.  For convenience, let $c_0=0$.

In view of the kernel method for the bivariate case, we locate the maximum $z^{max}$ of $z$ on $H(x,y,z)=0$. In order to do it, taking  the derivative with respect to $x$ at the both side of \eqref{s3-1} yields
\beqnn
-x\sigma_1^2-z\frac{\partial z}{\partial x} \sigma_3^2+(c_3-\lambda_3-c_{2})\frac{\partial z}{\partial x}+(c_1-\lambda_1-c_{0})=0,
\eeqnn
i.e.,
\beqlb\label{s2-22}
\big( -z\sigma_3^2-\lambda_3-c_{2}+c_3\big)\frac{\partial z}{\partial x}+\big(-x\sigma_1^2-\lambda_1-c_{0}+c_1\big)=0.
\eeqlb
Let
\beqlb\label{s3-2}
\frac{\partial z}{\partial x}=0,
\eeqlb
and solve the system of  equations (\ref{s2-22}) and \eqref{s3-2}, we have
\beqlb\label{s2-26}
x_{z^{max}}=\frac{c_1-\lambda_1-c_{0}}{\sigma^2_1}.
\eeqlb
Similarly,  take the derivative with respect to $y$,
\beqlb\label{s3-2a}
\frac{\partial z}{\partial y}=0,
\eeqlb
to obtain
\beqlb\label{s2-26-a}
y_{z^{max}}=\frac{c_2-\lambda_2-c_{1}}{\sigma^2_2}.
\eeqlb
It is easy to check that at the point $\big(x_{z^{max}},y_{z^{max}}\big)$, $z$ attains the maximum value $z^{max}$. From (\ref{s2-26}) and \eqref{s2-26-a},  we can get that on the point $(x_{z^{max}},y_{z^{max}},z^{max})$, the coordinates $x$ and $y$ satisfy
\beqlb\label{s2-28}
x_{z^{max}}=k_1 y_{z^{max}},
\eeqlb
where \beqlb\label{s2-84}
k_1=\frac{(c_1-\lambda_1-c_0)\sigma_2^2}{\big(c_2-\lambda_2-c_1\big)\sigma_1^2}.
\eeqlb
\begin{rem}
Without loss of generality, we  assume that $k_1\neq 1$ in the rest of this paper.  For the special case $k_1=1$, the discussion can be carried out by using the same ideal which is much simpler than the general case due to the fact that when $k_1=1$, the term including $k_1-1$ in most equations will disappear.
\end{rem}

From the above arguments, we obtain the maximum $z^{max}$ on the plane $H(k_1y,y,z)=0$.
Now, we consider  the new equation:
\beqlb\label{s2-86}
 H(k_1y,y,z)=0.
\eeqlb
From \eqref{s2-6} and \eqref{s3-1}, we can easily know that \eqref{s2-86} defines an ellipse.
 Thus, for fixed $z$, there are two solutions to \eqref{s2-86} for $y$, which are given:
\beqlb\label{s2-87}Y_{max,0}(z)&&=\frac{(c_1-\lambda_1) k_1+(c_2-\lambda_2-c_1)-\sqrt{
\Delta(z)}}{(\sigma_1^2 k_1^2+\sigma_2^2)},\;\textrm{and}
\\
\label{s2-88}
Y_{max,1}(z)&&=\frac{(c_1-\lambda_1) k_1+(c_2-\lambda_2-c_1)+\sqrt{
\Delta(z)}}{(\sigma_1^2 k_1^2+\sigma_2^2)},
\eeqlb
where
\beqlb\label{s2-89}
\Delta(z)=\Big(\big(c_1-\lambda_1 \big) k_1+c_2-\lambda_2-c_1\Big)^2+2\big(\sigma_1^2 k_1^2+\sigma_2^2\big)\big(-\frac{1}{2}\sigma_3^2z^2+(c_3-\lambda_3-c_2)z\big).
\eeqlb
Moreover, these two solution are distinct  except $\Delta(z)=0$.
We call a point $z$ a branch point if   $\Delta(z)=0$.   For branch points, we have the following property.
\begin{lem}\label{lem-10}\quad
\begin{itemize}
\item[(i)]$\Delta(z)$ has two real zeros, one of which is $z^{max}$, and the other is denoted by $z^{min}$. Moreover, they satisfy
\beqlb\label{s2-78}
z^{min}<0<z^{max}.
\eeqlb
\item[(ii)]
 $\Delta(z)>0$ in $(z^{min},\;z^{max})$ and  $\Delta(z)<0$ in $(-\infty,\;z^{min})\cup(z^{max},\;\infty)$.
\end{itemize}
\end{lem}
{\it Proof: } From \eqref{s2-89}, we obtain
\beqlb\label{s2-80}
\Delta(0)=\Big(\sum_{i=1}^{2} (c_i-\lambda_i-c_{i-1})k_i\Big)^2>0,
\eeqlb
where $k_2=1$. On the other hand,
\beqlb\label{s2-81}
\sum_{i=1}^{2}\sigma_i^2 k_i^2>0.
\eeqlb
From (\ref{s2-80}) and (\ref{s2-81}), we get \eqref{s2-78}.

By  properties of quadratic functions, we can get that (ii) holds.  The proof of the lemma is completed now.
\qed

In order to use the Tauberian-like Theorem below, we consider the analytic continuation of the moment generating functions in the complex plane $\mathbb{C}$.  The function  $\sqrt{\Delta(x)}$ plays an important role in the procedure of  the analytic continuation.  Hence, we first study its analytic continuation.
By Lemma \ref{lem-10}, $\sqrt{\Delta(x)}$ is well defined for $x\in[z^{min},\;z^{max}]$.   Moreover, it is a multi-valued function in the complex plane.  For convenience, in the sequel, $\sqrt{\Delta(x)}$ denotes the principle branch, that is $\Delta(x)=\Delta(\Re (x))$ for $x\in(z^{min},\;z^{max})$. In the follow, we continue $\sqrt{\Delta(x)}$ to the cut plane $\C\setminus\big\{(-\infty,\;z^{min}]\cup[z^{max},\infty)\big\}$.  In fact, we have
\begin{lem}\label{s3-lem1}
$\sqrt{\Delta(x)}$ is analytic in the cut plane $\C\setminus\big\{(-\infty,\;z^{min}]\cup[z^{max},\infty)\big\}$.
\end{lem}
The proof of Lemma \ref{s3-lem1} is standard. For example, see Dai and Miyazawa \cite{DM2011}, or Dai, Dawson and Zhao \cite{DDZ}.  However, for the completeness of the paper, we  provide a proof following the ideal used by   Dai and Miyazawa \cite{DM2011}.

\noindent{\it Proof:}  Since $z^{min}$ and $z^{max}$ are two zeros of $\Delta(x)=0$, we  have
\beqlb\label{s3-8}
\Delta(x)=\bigg(\sum_{i=1}^{2}\sigma_i^2 k_i^2\bigg)\sigma_3^2(x-z^{min})(z^{max}-x).
\eeqlb
Next, we rewrite \eqref{s3-8} in the polar form. Let $\omega_{min}(x)$ and $\omega_{max}(x)$ denote the principal arguments of $x-z^{min}$ and $z^{max}-x$, respectively.  Therefore,
\beqlb\label{s3-16}
\omega_{min}(x),\omega_{max}(x)\in(-\pi,\;\pi].
\eeqlb
Hence, \eqref{s3-8} can be rewritten as
\beqlb\label{s3-9}
\Delta(x)=\bigg[\sum_{i=1}^{2}\sigma_i^2 k_i^2\bigg]\sigma_3^2|x-z^{min}||z^{max}-x|\exp\{i(\omega_{min}(x)+\omega_{max}(x))\}.
\eeqlb
Moreover for $\omega_{min},\omega_{max}\in(-\pi,\;\pi]$, the functions  $x-z^{min}$ and $z^{max}-x$ are analytic.
Since $\sqrt{\Delta(x)}$ is the principle part, we have
\beqlb\label{s3-27}
\sqrt{\Delta(x)}=\sqrt{\big(\sum_{i=1}^{2}\sigma_i^2 k_i^2\big)\sigma_3^2|x-z^{min}||z^{max}-x|}\exp\Big[i\frac{\omega_{min}(x)+\omega_{max}(x)}{2}\Big].
\eeqlb
Thus, from \eqref{s3-16} and \eqref{s3-27}, one can easily get that $\sqrt{\Delta(x)}$ is analytic in the cut plane.\qed
\begin{cor}\label{s3-cor1}
Both $Y_{max,0}(x)$ and $Y_{max,1}(x)$ are analytic in the cut plane $\C\setminus\big\{(-\infty,\;z^{min}]\cup[z^{max},\infty)\big\}$.
\end{cor}

Symmetrically, we can treat the kernel equation \eqref{s2-86} as a quadratic function in $z$, and obtain parallel results to those in Lemmas \ref{lem-10} and \ref{s3-lem1}, and  Corollary \ref{s3-cor1}. We list them below.
Before  stating them, we first introduce the following notation. Define
\beqlb\label{s5-3}
\bar{\Delta}(y)=(c_3-\lambda_3-\lambda_{2})^2+2\sigma_3^2\Big(y\sum_{i=1}^{2}(c_i-\lambda_i-c_{i-1})k_i-\frac{1}{2}y^2\sum_{i=1}^{2}\sigma_i^2k_i^2\Big).
\eeqlb
For fixed $y$, there are two solutions to \eqref{s2-86}, which are given by
\beqlb\label{s5-1a}
&&Z_{max,1}(y)=\frac{(c_3-\lambda_3-c_{2})+\sqrt{\bar{\Delta}(y)}}{\sigma_3^2},\;\textrm{and}
\\
\label{s5-2a}
&&Z_{max,0}(y)=\frac{(c_3-\lambda_3-c_{2})-\sqrt{\bar{\Delta}(y)}}{\sigma_3^2}.
\eeqlb

Similar to Lemmas \ref{lem-10} and \ref{s3-lem1}, and  Corollary \ref{s3-cor1}, we have:
\begin{lem}\label{lem-10a}\quad
\begin{description}
\item [(i)]$\bar{\Delta}(y)$ has two real zeros, denoted by $y^{min}$ and $y^{max}$, respectively, satisfying
\beqlb\label{s2-78a}
y^{min}<0<y^{max}.
\eeqlb
\item[(ii)] $\bar{\Delta}(y)>0$ in $(y^{min},\;y^{max})$ and  $\bar{\Delta}(y)<0$ in $(-\infty,\;y^{min})\cup(y^{max},\;\infty)$.
    \item[(iii)] $Z_{max,0}(y)$ are analytic in the cut plane $\mathbb{C}\setminus\{(-\infty,y^{min}]\cup[y^{max},\infty)\}$.
    \end{description}
\end{lem}

In order to get the analytic continuation of the moment generating functions, we need  some technical lemmas. Before we introduce these lemmas, we first present an important notation. Define
\beqnn
\mathbb{G}_{\delta}(z_0)=\big\{z\in\C:z\neq z_0, |arg(z-z_0)|>\delta\big\}\cap\big\{z:|z|<z_0+\epsilon \big\},
\eeqnn
where $arg(z)\in(-\pi,\;\pi)$, $\delta\geq0$ and $\epsilon>0$.

For the function $Y_{max,0}(x)$, we have the following properties.
\begin{lem}\label{lem-9} For $Y_{max,0}(x)$, we have
\begin{itemize}
\item[(i)] $\Re~(Y_{max,0}(z))\leq Y_{max,0}\big(\Re (z)\big)$ for $\Re (z)\in(z^{min},\;z^{max})$.
\item[(ii)] $\Re (Y_{max,0}(z))<y^{max}$ for $z\in\mathbb{G}_{\delta_0}(z^{max})\cap\{z\in\C:\;z^{min}<\Re(z)\}$ with some $\delta_0\in[0,\;\frac{\pi}{2})$.
\end{itemize}
\end{lem}
{\it Proof:}~ It follows from \eqref{s2-87} and \eqref{s3-8} that
\beqlb\label{s3-10}
&&\Re~(Y_{max,0}(z))-Y_{max,0}\big(\Re(z)\big)\nonumber
\\&&\quad=\sqrt{\sigma_3^2\sum_{i=1}^{2}\sigma_i^2k_i^2}\bigg[\sqrt{(\Re(z)-z^{min})(z^{max}-\Re(z))}-\Re\Big(\sqrt{ (z-z^{min})(z^{max}-z)}\Big)\bigg].
\eeqlb
By \eqref{s3-10}, in order to prove  case (i), we only need to show
\beqlb\label{s3-11}
\sqrt{(\Re(z)-z^{min})(z^{max}-\Re(z))}-\Re\Big(\sqrt{ (z-z^{min})(z^{max}-z)}\Big)\leq 0.
\eeqlb
We also note that  $(\Re(z)-z^{min})$ and $(z^{max}-\Re(z))$ are real parts of $(z-z^{min})$ and $(z^{max}-z)$, respectively, since $z^{min}$ and  $z^{max}$  are real.  Therefore,
\beqnn
(\Re(z)-z^{min})&&=|z-z^{min}|\cos\big(\omega_{min}(z)\big),
\\
\big(z^{max}-\Re(z)\big)&&=|z^{max}- z|\cos\big(\omega_{max}(z)\big).
\eeqnn
So,
\beqlb\label{s3-20}
\sqrt{(\Re(z)-z^{min})(z^{max}-\Re(z))}&&=\sqrt{|z-z^{min}||z^{max}- z|}
\Big(\cos\big(\omega_{min}(z)\big)\cos\big(\omega_{max}(z)\big)\Big)^{\frac{1}{2}}.\qquad
\eeqlb
Similarly, we have
\beqnn
(z-z^{min})(z^{max}-z)&&=|z-z^{min}||z^{max}- z|\exp\Big\{i(\omega_{max}(z)+\omega_{min}(z))\Big\}.
\eeqnn
Thus,
\beqnn
\sqrt{(z-z^{min})(z^{max}-z)}&&=\sqrt{|z-z^{min}||z^{max}- z|}\nonumber
\exp\{i\frac{\omega_{max}(z)+\omega_{min}(z)}{2}\}.
\eeqnn
Hence,
\beqlb\label{s3-14}
\Re\Big( \sqrt{(z-z^{min})(z^{max}-z)}\Big)&&=\sqrt{|z-z^{min}||z^{max}- z|}\cos \big(\frac{\omega_{max}(z)+\omega_{min}(z)}{2}\big).
\eeqlb
Since for $\Re(z)\in(z^{min}\;z^{max})$,
\beqlb\label{s3-21}
&&\omega_{max}(z)\in(-\frac{\pi}{2},\;0),
\\\label{s3-22}
&&\omega_{min}(z)\in(0,\;\frac{\pi}{2}),
\\
\label{s3-15}
&&\frac{\omega_{max}(z)+\omega_{min}(z)}{2}\in(-\frac{\pi}{4},\;\frac{\pi}{4}).
\eeqlb
From \eqref{s3-20} to \eqref{s3-15}, in order to prove \eqref{s3-11}, we only need to prove
\beqlb\label{s3-17}
\Big(\cos\big(\omega_{min}(z)\big)\cos\big(\omega_{max}(z)\big)\Big)^{\frac{1}{2}}\leq \cos\Big(\frac{\omega_{max}(z)+\omega_{min}(z)}{2}\Big),\qquad
\eeqlb
which directly follows from Dai and Miyazawa \cite{DM2011}.

Next, we prove  case (ii). We first assume that $z^{min}<\Re(z)<z^{max}$. From \eqref{s2-87} and Lemma \ref{lem-10}, we have
\beqlb\label{s3-23}
y^{m}:=Y_{max,0}(z^{max})=\frac{\sum_{i=1}^{2}(\lambda_i-c_i-c_{i-1}) k_i}{\sum_{i=1}^{2}\sigma_i^2 k_i^2},
\eeqlb
since $\Delta(z^{max})=0$.  From \eqref{s3-23} and  case (i), in order to prove the case (ii), we only need to show
\beqlb\label{s3-12}
\frac{\sum_{i=1}^{2}(c_i-\lambda_i-c_{i-1})k_i}{\sum_{i=1}^{2}\sigma_i^2 k_i^2}\leq y^{max}.
\eeqlb
On the other hand, it follows from  \eqref{s5-3} and Lemma \ref{lem-10a} that
\beqlb\label{s3-13}
2 \frac{\sum_{i=1}^{2}(\lambda_i-c_i-c_{i-1}) k_i}{\sum_{i=1}^{2}\sigma_i^2 k_i^2}\leq y^{max}.
\eeqlb
Hence, \eqref{s3-12} follows from \eqref{s3-13}.

Finally, we assume that $\Re(z)\geq z^{max}$.   As $\delta\to \frac{\pi}{2}$, we have that
\beqlb\label{s3-3}
\Re(z)\to z^{max}.
\eeqlb
It follows from Lemma \ref{lem-10a}, \eqref{s3-13} and \eqref{s3-3} that we can find $\delta_0\in(0,\;\frac{\pi}{2})$ such that case (ii) holds.
The proof of the lemma is completed.
 \qed
\subsection{Analytic Continuation}
The analytic continuation of  the moment generating function $\phi_2(0,0,z)$ plays an important role in our analysis, which is the focus in  this subsection.  In order to carry out this, we need the following technical lemma.
\begin{lem}\label{lem-a-5}
\label{lem-11}
For the  moment generating functions $\phi_i(\cdot)$,$i=1,2,3$, we have
\begin{itemize}
\item[(i)] $\phi_1(0,y,z)$ is finite on some region $\{(y,z):y<\epsilon,z<\epsilon\}$ with $\epsilon>0$.
\item[(ii)] $\phi_2(0,0,z)$ is finite on some region $\{z: z<\epsilon\}$ with $\epsilon>0$.
\item[(iii)] $\phi_2(x,0,z)$ is finite on some region $\{z: z<\epsilon,x<\epsilon\}$ with $\epsilon>0$.
\item[(iv)] $\phi_3(0,y,0)$ is finite on some region $\{y: y<\epsilon\}$ with $\epsilon>0$.
\end{itemize}
\end{lem}
{\it Proof:} We first prove  case (i).  In order to prove it, we first prove
\beqlb\label{s4-6}
\E\Big[\int_0^1 e^{y L_2(u)}dY_1(u)\Big]<\infty
\eeqlb
for some $y>0$, and
\beqlb\label{s4-8}
\E\Big[\int_0^1 e^{zL_3(u)}d Y_1(u)\Big]<\infty
\eeqlb
for some $z>0$.

 In fact,
\beqlb\label{s4-3}
 \E\Big[\int_0^1 e^{y L_2(u)}dY_1(u)\Big]=\phi_1(0,y,0).
\eeqlb
which suggests that  we may  restrict our analysis to the two-dimensional tandem queue $\{( L_1(t), L_2(t))'\}$ with the two nodes 1 and 2. We  note that $( L_1(t), L_2(t))'$ is not affected by $L_3(t)$.  Using the same method as in Dai, Dawson and Zhao \cite{DDZ}, we can easily get \eqref{s4-6}.

Next, we prove \eqref{s4-8}. Since $Y_1$ is a regulator,
\beqlb\label{s4-8a}
\E\Big[\int_0^1 e^{zL_3(u)}d Y_1(u)\Big]=\E\Big[\int_0^1 e^{xL_1(t)+0L_2(t)+zL_3(u)}d Y_1(u)\Big]=\phi_1(x,0,z).
\eeqlb
By \eqref{s2-69} and \eqref{s4-8a}, we get that the left-hand side of \eqref{s4-8} satisfies
\beqlb\label{s4-9}
H(x,0,z)\phi(x,0,z)=x\phi_1(0,0,z)-z\phi_2(x,0,z)+z\phi_3(x,0,0).
\eeqlb

Next, we study this system on the plane $y=0$. We first consider the  ellipse defined by
\beqlb\label{s4-10}
H(x,0,z)=0.
\eeqlb
For the point $(x,z)$ on this ellipse, we have
\beqlb\label{s4-11}
x\phi_1(0,0,z)-z\phi_2(x,0,z)+z\phi_3(x,0,0)=0.
\eeqlb
For fixed $x$, we can find two solutions to \eqref{s4-10} for $z$. Denote one  of these two solutions by
\beqlb\label{s4-12}
Z_0(x)=\frac{(c_3-\lambda_3-c_2)-\sqrt{(c_3-\lambda_3-c_2)^2+2\sigma_3^2(-\frac{1}{2}\sigma_1^2x^2+(c_1-\lambda_1)x)}}{\sigma^2_3}.
\eeqlb

Using the same method as in the proof of Lemma \ref{lem-10}, we can get that $Z_0(x)$ is well-defined between $[x^{min},\;x^{max}]$ with $x^{min}<0<x^{max}$ and
\beqnn
\Delta_1(x^{min})=\Delta_1(x^{max})=0,
\eeqnn
where
$$
\Delta_1(x)=(c_3-\lambda_3-c_2)^2+2\sigma_3^2(-\frac{1}{2}\sigma_1^2x^2+(c_1-\lambda_1)x).
$$
Hence, from \eqref{s4-11} and \eqref{s4-12}, we have
\beqlb\label{s4-13}
x\phi_1(0,0,Z_0(x))-Z_0(x)\phi_2(x,0,Z_0(x))+Z_0(x)\phi_3(x,0,0)=0,
\eeqlb
that is,
\beqlb\label{s4-14}
Z_0(x)\phi_3(x,0,0)=Z_0(x)\phi_2\big(x,0,Z_0(x)\big)-x\phi_1\big(0,0,Z_0(x)\big).
\eeqlb
Hence,   $Z_0(x)\phi_3(x,0,0)$ is finite if and only if the right-hand side of \eqref{s4-14} is finite.
On the other hand, from \eqref{s4-12}, we obtain that for  $x\in[x^{min},\;0)$,
\beqlb\label{s4-15}
z=Z_0(x)>0,
\eeqlb
and
\beqlb\label{s4-16}
\phi_3(x,0,0)<\infty.
\eeqlb
From \eqref{s4-15} and \eqref{s4-16}, we obtain that
\beqlb\label{s4-17}
\phi_1(0, 0,Z_0(x))<\infty,
\eeqlb
since $ \infty>Z_0(x)\phi_2\big(x,0,Z_0(x)\big)\geq 0$ and $-x\phi_1\big(0,0,Z_0(x)\big)\geq 0$. Therefore   \eqref{s4-8} holds. Finally, we have
\beqlb\label{s4-18}
\phi_1(0,y,z)&&=\E\Big[\int_0^1 e^{yL_2(u)+zL_3(u)}dY_1(u)\Big]\nonumber
\\&&\leq \frac{1}{2}\Big(\E\Big[ \int_0^1  e^{2yL_2(u)}d Y_1(u)\Big]+ \E\Big[\int_0^1  e^{2zL_3(u)}d Y_1(u)\Big]\Big).
\eeqlb
Combing \eqref{s4-6},  \eqref{s4-8} and \eqref{s4-18}, we get that for some $y>0$ and $z>0$
\beqnn
\phi_1(0,y,z)<\infty.
\eeqnn

Next, we prove  case (ii).   Since
\beqlb
\phi_2(0,0,z)=\E\Big[\int_0^1 e^{0L_1(u)+zL_3(u)}dY_2(u)\Big],
\eeqlb
we can  consider the problem on the plane $x=0$. It follows from \eqref{s2-69} that
\beqlb\label{s4-19}
H(0,y,z)\phi(0,y,z)=-y\phi_1(0,y,z)+(y-z)\phi_2(0,0,z)+z\phi_3(0,y,0).
\eeqlb
Then,
  \beqlb\label{s3-26}
H(0,y,z)=0
\eeqlb defines an ellipse.  For every fixed $y$,
define
\beqlb\label{s2-44-a}
\bar{Z}_0(y)=\frac{(c_3-\lambda_3-c_2)-\sqrt{(c_3-\lambda_3-c_2)^2+2\sigma_3^2(-\frac{1}{2}\sigma_2^2y^2+(c_2-\lambda_2-c_1)y)}}{\sigma_3^2}.
\eeqlb
Then, \eqref{s2-44-a} is a solution to  equation \eqref{s3-26}. Similar to Lemma \ref{lem-10}, $\bar{Z}_0(y)$  is well-defined on some region $[a,\;b]$ with $a<0$ and $b>0$. It follows from \eqref{s4-19} and \eqref{s2-44-a} that
\beqlb\label{s2-98}
\big(y-\bar{Z}_0(y)\big)\phi_2(0,0,z)=y\phi_1\big(0,y,\bar{Z}_0(y)\big)-\bar{Z}_0(y)\phi_3(0,y,0).
\eeqlb
Furthermore, from \eqref{s2-44-a}, we obtain that for $y\in\{y:a<y<0\}$
\beqlb\label{s4-20}
\bar{Z}_0(y)>0.
\eeqlb
Hence, by  case (i) and \eqref{s4-20}, we can choose $y<0$ such that $z=\bar{Z}_0(y)>0$ and,
\beqlb
\label{s4-22}
\phi_1\big(0,y,\bar{Z}_0(y)\big)<\infty.
\eeqlb
It is also worthy noting that  for $y<0$,
\beqlb
\label{s4-21}
\phi_3(0,y,0)<\infty.
\eeqlb
Case (ii) now follows from  \eqref{s2-98} to \eqref{s4-21}.

Finally, we can show  cases (iii) and (iv)to complete the proof of the lemma.
\qed

For the continuation of the function  $\phi_2(0,0,z)$,  we need another technical tool.

\begin{lem}\label{s4-lem1}
Let $f(x_1,x_2)$ be a probability density function on $\R^2_+$.
For a real variable $\lambda$, define $\tilde{G}(\lambda)=\int_{\R_+^2} e^{g(\lambda)x_1+\lambda x_2 }f(x_1,x_2)dx$ with $g(\lambda)$ being a bounded and continuously differential real function, and
\beqlb\label{4-1}
\tau_{\tilde{G}}=\sup\{\lambda\geq 0: \tilde{G}(\lambda)<\infty\}.
\eeqlb
Then, the complex variable function $\tilde{G}(z)$ is analytic on $\{z\in \mathbb{C}: \Re (z)<\tau_{\tilde{G}}\}$.
\end{lem}
{\it Proof:} We use the Vitali's Theorem to prove it.
In fact, we have
\beqlb\label{s4-54}
\int_{\R_+^2} e^{g(\lambda)x_1+\lambda x_2 }f(x_1,x_2)dx_1dx_2=\int_0^\infty e^{g(\lambda)x_1}dx_1 \Big[\int_0^\infty e^{\lambda x_2}f(x_1,x_2)dx_2\Big].
\eeqlb
For convenience, define
\beqlb\label{s4-55}
F(\lambda,x_1)=\int_0^\infty e^{\lambda x_2}f(x_1,x_2)dx_2.
\eeqlb
Since $f(x_1,x_2)$ is a density function,
we can get that $F(\lambda,x_1)$ is analytic on the region $\{z\in \mathbb{C}: \Re (z)<\tau_{\tilde{G}}\}$ for any $x_1\in\R_+$. Let
\beqlb\label{s4-62}
\tilde{F}(\lambda,x_1)=e^{g(\lambda)x_1}F(\lambda,x_1).
\eeqlb
Now, it is obvious that $\tilde{F}(\lambda,x_1)$ satisfies the conditions of the Vitali's Theorem (see, for example, Markushevich \cite{M1977}) on the region $\{z\in \mathbb{C}: \Re (z)<\tau_{\tilde{G}}\}$. Then, the lemma holds. \qed

\begin{rem} From Lemma \ref{s4-lem1},
\begin{itemize}
\item[(i)] The convergence parameter $\tau_G$ is unique;
\item[(ii)] If $\alpha\in\R_+$ is a singularity of $G$, then we must have  $G(x)=\infty$ for $x\in(\alpha,\infty)$. However, $G(\alpha)$ may be either finite or infinite.
\end{itemize}
\end{rem}
\begin{rem}\label{rem-3}
It follows from Lemmas \ref{lem-11} and \ref{s4-lem1} that
\begin{itemize}
\item[(i)] $\phi_2(0,0,z)$ is finite on some region $\{\ z: \Re z<\epsilon\}$ with $\epsilon>0$, which implies that the convergence parameter $\tau_{\phi_2(0,0,z)}$ is greater than 0.
\item[(ii)] $\phi_3(0,y,0)$ is finite on some region $\{y: \Re y<\epsilon\}$ with $\epsilon>0$, which implies that the convergence parameter $\tau_{\phi_3(0,y,0)}$ is greater than 0.
\end{itemize}
\end{rem}

The next lemma enables us to express $\phi_{2}(0,0,z)$ in terms of the other moment generating functions.

\begin{lem}\label{lem-a-4}
$\phi_2(0,0,z)$  can be analytically continued to the region $z\in\{z: \Re z<\epsilon\}$ with $\epsilon>0$, and
\beqlb\label{s2-74}
\phi_2(0,0,z)=&&\frac{\phi_2(k_1 Y_{max,0}(z),0,z)}{1-k_1}\nonumber
\\&&\quad+z \frac{\phi_3(k_1 Y_{max,0}(z),Y_{max,0}(z),0)}{(1-k_1)( Y_{max,0}(z)-z)}-z \frac{\phi_3(0,0, Y_{max,0}(z))}{( Y_{max,0}(z)-z)}.
\eeqlb

\end{lem}

\noindent{\it Proof:}
From Corollary \ref{s3-cor1} and \eqref{s2-69}, we get that
\beqlb\label{s2-90}
&&\phi_1\big(k_1Y_{max,0}(z),Y_{max,0}(z),z\big)=\nonumber
\\&&\qquad-\frac{\phi_2\big(k_1Y_{max,0}(z),0,z\big)(Y_{max,0}(z)-z)+z\phi_3\big(k_1Y_{max,0}(z),Y_{max,0}(z),0\big)}{(k_1-1)Y_{max,0}(z)}.
\eeqlb
On the other hand, equation \eqref{s3-26}  defines an ellipse. For fixed $z$, there are two solutions to \eqref{s3-26} for $y$.  Define
\beqlb\label{s4-23}
Y_{0}(z)=\frac{(c_2-\lambda_2-c_1)-\sqrt{\Delta_{H(0,y,z)}(z)}}{\sigma_2^2},
\eeqlb
where  $$\Delta_{H(0,y,z)}(z)=(c_2-\lambda_2-c_1)^2+2\sigma_2^2\big((c_3-\lambda_3-c_2)z-\frac{1}{2}\sigma_3^2z^2\big).$$
Using the same method as in  the proof of Lemma \ref{s3-lem1}, we can get that $Y_{0}(z)$ is analytic in the cut plane $\mathbb{C}\setminus\{(-\infty,\bar{z}^{min}]\cup[\bar{z}^{max},\;\infty)\}$, where
\beqnn
\Delta_{H(0,y,z)}(\bar{z}^{min})=\Delta_{H(0,y,z)}(\bar{z}^{max})=0
\eeqnn
with
\beqnn
\bar{z}^{min}<0<\bar{z}^{max}.
\eeqnn
By (\ref{s4-19}) and \eqref{s4-23}, We can find a region such that
\beqlb\label{s2-41-a}
\phi_1(0,Y_{0}(z),z)=\frac{\phi_2(0,0,z)(Y_{0}(z)-z)+z\phi_3(0,Y_{0}(z),0)}{Y_{0}(z)}.
\eeqlb
Next, we  study the relationship between $Y_{0}(z)$ and $Y_{max,0}(z)$ for  $z>0$.
We  note that  both the two ellipses defined by \eqref{s2-86} and \eqref{s3-26}, respectively, pass the origin $(0,0)$ and
\beqlb\label{s2-43-a}
H(k_1y,y,z)&&=H(0,y,z)-\frac{1}{2}k_1^2\sigma_1^2y^2+(c_1-\lambda_1) k_1y\nonumber
\\&&=H(0,y,z)+G(y),
\eeqlb
where $$G(y)=-\frac{1}{2}k_1^2\sigma_1^2y^2+(c_1-\lambda_1) k_1y.$$

We should note that
\beqlb\label{s2-56-a}
G(y)>0\;\textrm{if and only if }\;y\in\Big[0, \frac{2(c_1-\lambda_1)}{\sigma_2^2k_1}\Big].
\eeqlb
From \eqref{s4-23}, we obtain that for $z\in[0,\;\frac{2(c_3-\lambda_3-c_2)}{\sigma_3^2}]$
\beqlb\label{s4-2}
Y_{0}(z)\leq 0.
\eeqlb
From \eqref{s2-43-a} and \eqref{s4-2}, we get that for $z\in(0,\;\frac{2(c_3-\lambda_3-c_2)}{\sigma_3^2})$
\beqlb\label{s2-46-a}
H(k_1Y_{0}(z),Y_{0}(z),z)<0.
\eeqlb
On the other hand, from \eqref{s2-87} and \eqref{s4-23}, we have, for $0<z<\frac{2(c_3-\lambda_3-c_2)}{\sigma_3^2},$
\beqnn
Y_{0}(z)&&<0\;\textrm{and}\;
Y_{max,0}(z)<0.
\eeqnn
Thus,
\beqlb\label{s2-59}
Y_{max,0}(z)>Y_{0}(z),
\eeqlb
for $z\in\Big(0, \frac{2(c_3-\lambda_3-c_2)}{\sigma_3^2}\Big)$.

It follows from  Lemmas \ref{lem-11}, \eqref{s2-90}, and \eqref{s2-41-a} that
\beqnn\label{s2-79}
\phi_1(0,Y_{max,0}(z),z)&&=\frac{\phi_2(0,0,z)(Y_{max,0}(z)-z)+z\phi_3(0,Y_{max,0}(z),0)}{Y_{max,0}(z)}\nonumber
\\&&=-\frac{\phi_2\big(k_1Y_{max,0}(z),0,z\big)(Y_{max,0}(z)-z)+z\phi_3(k_1Y_{max,0}(z),Y_{max,0}(z),0)}{(k_1-1)Y_{max,0}(z)},\qquad
\eeqnn
where we use the principle  of  analytic continuation of several complex variables functions (see, for example, Narasimhan \cite{N1964}).
Therefore
\beqlb\label{s2-96}
\phi_2(0,0,z)=&&\frac{\phi_2(k_1Y_{max,0}(z), 0, z)}{1-k_1}\nonumber
\\&&+z \frac{\phi_3(k_1Y_{max,0}(z),Y_{max,0}(z),0)}{(1-k_1)(Y_{max,0}(z)-z)}-z \frac{\phi_3(0,Y_{max,0}(z),0)}{(Y_{max,0}(z)-z)},
\eeqlb
 for $\Re~z<\epsilon$ with some $\epsilon>0$. The proof is completed.
\qed

\section{Tail Asymptotic and Singularity Analysis }

In order to reach our goal, we need to study  tail behaviors of $\phi_2(0,0,z)$ around the dominant singularities.   From Lemma \ref{s4-lem1}, there exists only one dominant singularity.  We denote it by $z_{dom}$. Next, we  characterize the dominant singularity $z_{dom}$ of $\phi_2(0,0,z)$. For convenience, let
\beqnn
F(y)&&= \frac{\phi_3(k_1 y,y,0)}{(1-k_1)}- \phi_3(0,y,0)\;\textrm{and}\; D(y,z)=\phi_2(0,y,z)-\frac{\phi_2(k_1y,y,z)}{1-k_1}.
\eeqnn
Moreover, let
\beqlb\label{s4-d-4}
G(z)=D(Y_{max,0}(z),z),\;\textrm{and}\;\bar{G}(y)=D(y,Z_{max,0}(y)).
\eeqlb
From Lemma \ref{lem-a-4}, we have:
\begin{lem}\label{lem-2} $F(y)$  can be analytically continued to a region $\{y:\Re(y)<\epsilon\}$ with $\epsilon>0$, and
\beqlb\label{s2-61a}
F(y)=\frac{Z_{max,0}(y)-y}{Z_{max,0}(y)}
\bar{G}(y).
\eeqlb
\end{lem}

 We introduce the following notation.
\beqnn
\hat{\phi}_2(0,0,z)&&=\phi_2(k_1Y_{max,0}(z),0,z).
\eeqnn
Next, we first study the relationship between the convergence parameters of $\phi_2(0,0,z)$, $\hat{\phi}_{2}(0,0,z)$ and $G(z)$. In fact, we have:

\begin{lem}\label{lem-d-2} For the convergence parameters  $\tau_{\phi_2}$, $\tau_{\hat{\phi}_2}$ and $\tau_G$ of $\phi_2(0,0,z)$, $\hat{\phi}_2(0,0,z)$ and $G(z)$, respectively,  we have
\beqlb\label{s4-d-1}\tau_G=\tau_{\phi_2}=\tau_{\hat{\phi}_2}.\eeqlb
\end{lem}
{\it Proof:} We first show
\beqlb\label{s4-d-2}\tau_{\phi_2}=\tau_{\hat{\phi}_2}.  \eeqlb
By Lemma \ref{lem-a-5}, we just need to focus on  $z>0$. By \eqref{s2-10}, we get that if $Y_{max,0}(z)\geq 0$, then
\beqlb\label{s4-25}\tau_{\phi_2}\geq\tau_{\hat{\phi}_2};\eeqlb if $Y_{max,0}(z)< 0$, then \beqlb\label{s4-42}\tau_{\phi_2}\leq \tau_{\hat{\phi}_2},\eeqlb  since $L_1(u)\geq0, L_2(u)\geq 0$, $z>0$ and $k_1>0$.

In order to prove \eqref{s4-d-2}, we first locate the dominant singularity $z_{dom}$.
From \eqref{s2-74}, we have
\beqlb\label{s4-67}
&&\phi_2(0,0,z)-\frac{\phi_2(k_1 Y_{max,0}(z), 0, z)}{1-k_1}\\\nonumber
\\&&\qquad\qquad=z \frac{\phi_3(k_1 Y_{max,0}(z),Y_{max,0}(z)(z),0)}{(1-k_1)(Y_{max,0}(z)-z)}-z \frac{\phi_3(0,Y_{max,0}(z),0)}{(Y_{max,0}(z)-z)}.
\eeqlb
We  observe from \eqref{s2-87} and \eqref{s2-89} that
\beqlb\label{s4-68}
Y_{max,0}\Big(\frac{2(c_3-\lambda_3-c_2)}{\sigma_3^2}\Big)=0.
\eeqlb
From \eqref{s2-85}, \eqref{s4-67} and \eqref{s4-68}, we get
\beqlb\label{s4-69}
\phi_2(0,0,\frac{2(c_3-\lambda_3-c_2)}{\sigma_3^2})=\phi_3(0,0,0)=c_3-\sum_{i=1}^3\lambda_i.
\eeqlb
Hence, from Lemma \ref{s4-lem1} and \eqref{s4-69}, we get
\beqlb\label{s4-97}
\tau_{\phi_2}> \frac{2(c_3-\lambda_3-c_2)}{\sigma_3^2}.
\eeqlb

For $z\in(\frac{2(c_3-\lambda_3-c_2)}{\sigma_3^2},\;z^{max})$, one can easily get that
\beqlb\label{s4-30}
Y_{max,0}(z)>0.
\eeqlb
Therefore,
\beqlb\label{s4-33}
\tau_{\phi_2}\geq \tau_{\hat{\phi}_2}.
\eeqlb
However, from \eqref{s2-96}, we must have \eqref{s4-d-2}.

Next, we prove \eqref{s4-d-1}. From  \eqref{s2-74} and \eqref{s4-d-2}, it is obvious that
\beqlb\label{s4-1}
\tau_{\phi_2}\leq z_{max}.
\eeqlb

If $\tau_{\phi_2}=z^{max}$, then, from \eqref{s2-96},
it must be the dominant singularity of $G(z)$.   Next, we assume $\tau_{\phi_2}\in(0,\;z^{max})$.  From \eqref{s2-10}, \eqref{s4-97} and \eqref{s4-30}, we have for $\frac{2(c_3-\lambda_3-c_2)}{\sigma_3^2}<z<z^{max}$
\beqlb\label{s4-72}
\hat{\phi}_2(0,0,z)-\phi_2\big(0,Y_{max,0}(z),z\big)&&=
\phi_2(k_1Y_{max,0}(z),0,z)-\phi_2(0,0,z)\nonumber
\\&&=\int_{\R_+^3}\Big(e^{k_1Y_{max,0}(z)x_{1}}-1\Big)e^{z x_{3}}V_{2}(dx)
\nonumber
\\&&
\geq \int_{\R_+^3}\Big(e^{k_1Y_{max,0}(z)x_{1}}-1\Big)V_{2}(dx)>0.
\eeqlb
It is worth  noting that,  from Lemma \ref{s4-lem1}, \eqref{s4-d-2},  and \eqref{s4-72}, we have
 \beqlb\label{s4-d-5}
 \lim_{z\to\tau_{\phi_2}}\phi_2(0,0,z)=\lim_{z\to\tau_{\phi_2}}\hat{\phi}_2(0,0,z)=\infty.
 \eeqlb
 If $\tau_{\phi_2}$ is not the dominant singularity of $G(z)$, then  $G(z)$ is analytic around $\tau_{\phi_2}$. So, $G(z)$ is bounded in a neighbourhood of $\alpha$.  On the other hand, from \eqref{s4-d-4}
\beqnn
\hat{\phi}_2(0,0,z)=(1-k_1)\phi_2(0,Y_{max,0}(z),z)-(1-k_1)G(z).
\eeqnn
Hence,
\beqlb\label{s4-d-7}
\hat{\phi}_2(0,0,z)-\phi_2(0,Y_{max,0}(z),z)=-k_1\phi_2(0,Y_{max,0}(z),z)-(1-k_1)G(z).
\eeqlb
From the  maximum modulus principle, Lemma \ref{s4-lem1} and \eqref{s4-d-7}, we obtain that for some region $\{z:0<|z-\alpha|\leq \epsilon\}$,
\beqlb\label{s4-d-8}
\hat{\phi}_2(0,0,z)-\phi_2(0,Y_{max,0}(z),z)<0,
\eeqlb
since $k_1>0$. It is obvious that \eqref{s4-72} contradicts to \eqref{s4-d-8}.  Hence the lemma holds. \qed
\begin{rem}
From the proof of Lemma \ref{lem-d-2}, we have the following important fact
\beqlb\label{s4-28}
\tau_{\phi_2}\leq z_{max}.
\eeqlb
\end{rem}

Next, we study the convergence parameter $\tau_G$.

\begin{lem}\label{s4-lem3}
If $z^{G}\in(0,\;z^{max}]$ is the dominant singularity of $G(z)$, then $\bar{G}(y)$ is analytic at the point $y^{0}:=Z_{max,0}(z^G)$.
\end{lem}
{\it Proof:} From \eqref{s5-2a}, we obtain the zero $y^*$ of $Z_{max,0}(y)$ is
\beqlb\label{s4-47}
y^*=2\frac{\sum_{i=1}^{2}(c_i-\lambda_i-c_{i-1})}{\sum_{i=1}^{2}\sigma_i^2k_i^2}.
\eeqlb
From  \eqref{s2-87} and Lemma \ref{lem-10}, we get
\beqlb\label{s4-48}
Y_{max,0}(z^{min})=Y_{max,0}(z^{max})=\frac{\sum_{i=1}^{2}(c_i-\lambda_i-c_{i-1})}{\sum_{i=1}^{2}\sigma_i^2k_i^2}:=\tilde{y}^m.
\eeqlb
Combing \eqref{s4-47} and \eqref{s4-48}, we obtain that
\beqlb\label{s4-d-13}
\tilde{y}^m<y^*.
\eeqlb
 It follows from \eqref{s5-3} and \eqref{s5-2a} that for $y\in(0,\;y^*)$
 \beqlb\label{s4-34}
 Z_{max,0}(y)<0.
 \eeqlb
 From \eqref{s2-87}, one can easily get that $Y_{max,0}(z)$ is increasing on $[\frac{2(c_3-\lambda_3-c_{2})}{\sigma_3^2},z^{max}]$.  Hence
 \beqlb\label{s4-24}
 y^0\leq \tilde{y}^m<y^*.
 \eeqlb
From \eqref{s4-34} and \eqref{s4-24}, we obtain that
 \beqlb\label{s4-36}
 Z_{max,0}(y^0)<0.
 \eeqlb
 Therefore $\phi_2(0,0,z)$  is analytic at the point $z^0:=Z_{max,0}(y^0)$.
 From \eqref{s4-d-4}, in order to prove the lemma, we only need to show that $\phi_2\big(k_1y,0,Z_{max,0}(y)\big)$ is analytic at $y^0$. From \eqref{s4-36}, we must have
 \beqlb\label{s4-38}
 Z_{max,1}(y^0)=z^G.
 \eeqlb
   It follows from Lemma \ref{lem-d-2} that $\phi_2(k_1Y_{max,0}(z),0,z)$ is analytic at $z^0$.  It follows from \eqref{s2-87} and \eqref{s4-48} that
\beqnn
Z_{max,0}(y^0)=z^0.
\eeqnn
From the above arguments, we can get that the lemma holds.
\qed

 The zero of $Y_{max,0}(z)-z$ is critical for us to prove Lemma \ref{lem-d-1}, Hence, we demonstrate how to evaluate it. Let $f(z)=\big(Y_{max,0}(z)-z\big)\big(Y_{max,1}(z)-z\big)$.  Then we have
\beqlb\label{s4-d-18}
f(z)=Y_{max,0}(z)Y_{max,1}(z)-z(Y_{max,1}(z)+Y_{max,0}(z))+z^2.
\eeqlb
It follows from  \eqref{s2-84}  and \eqref{s2-87} that
$$
-\frac{1}{2}(\sum_{i=1}^3\sigma_i^2k_i^2)f(z)=\Big(-\frac{1}{2}(\sum_{i=1}^{3}\sigma_i^2k_i^2)z+\sum_{i=1}^{3}(c_i-\lambda_i-c_{i-1})k_i\Big)z.
$$
Hence, the non-zero root of $Y_{max,0}(z)-z=0$  is
\beqlb\label{s2-97}
z^*=2\frac{\sum_{i=1}^3(c_i-\lambda_i-c_{i-1})k_i}{\sum_{i=1}^3\sigma_i^2k_i^2}.
\eeqlb

\begin{lem}\label{lem-d-1}
If $\tau_{G}\in(\frac{2(c_3-\lambda_3-c_2)}{\sigma_3^2},\;z^{max})$, then $\tau_{G}$ is the zero $z^*$ of $Y_{max,0}(z)-z$.
\end{lem}

\noindent{\it Proof:} From \eqref{s2-74}, we obtain that
\beqlb\label{s4-d-9}
G(z)=\frac{z}{Y_{max,0}(z)-z}F(Y_{max,0}(z)).
\eeqlb
Hence, in order to prove our result, we only need to show $F(Y_{max,0}(z))$ is analytic on $\{z:\Re(z)<z^*+\epsilon\}$ with small enough $\epsilon>0$. From  \eqref{s5-2a}, we have
\beqlb\label{s4-d-11}
Y_{max,0}(z^*)=z^*.
\eeqlb
Next, we show that
\beqlb\label{s4-d-12}
Y_{max,0}(z^*)\neq y^*.
\eeqlb
Since $Y_{max,0}(z)$ is increasing on $(\frac{2(c_3-\lambda_3-c_2)}{\sigma_3^2}, \; z^{max}]$, by \eqref{s4-d-13},
\beqlb\label{s4-d-14}
Y_{max,0}(z^*)<y^*.
\eeqlb
Finally, it follows from  Lemma \ref{s4-lem3} that $\bar{G}(y)$ is analytic at the point $Y_{max,0}(z^*)$.
From the above arguments and \eqref{s2-61a}, we have the lemma.\qed

From Lemma \ref{lem-d-1} and \eqref{s4-97}, we have:
\begin{lem}\label{lem-a1} If the convergence parameter $\tau_{\phi_2}$ is less than $z^{max}$, then
$$
\tau_{\phi_2}=2\frac{\sum_{i=1}^3(c_i-\lambda_i-c_{i-1})k_i}{\sum_{i=1}^3\sigma_i^2k_i^2}.
$$
\end{lem}

From Lemmas  \ref{lem-d-2} and \ref{lem-d-1}, we can get that $z_{dom}$ is either $z^*$ or $z^{max}$. In order to obtain  tail asymptotics for the marginal $L_3$, we need to study  asymptotic properties of the moment generating function $\phi_2(\cdot)$ at the point $z_{dom}$.  We first present asymptotic properties of $G(\cdot)$ at the point $z_{dom}$.
\begin{lem}\label{lem-d-3}
For the  function $G(z)$, we have
\begin{itemize}
\item[(i)] If $z_{dom}<z^{max}$, then $z_{dom}$ is a simple pole of $G(z)$, and
\beqlb\label{s4-d-15}
\lim_{z\to z_{dom}}(z_{dom}-z)G(z)=\frac{2(Y_{max,1}(z_{dom})-z_{dom})F(Y_{max,0}(z_{dom}))}{\sum_{i=1}^3\sigma_i^2k_i^2}.
\eeqlb
\item[(ii)]  If $z_{dom}=z_{max}<z^*$, then $z_{dom}$ is a branch point of $G(z)$. Moreover
\beqlb\label{s4-d-16}
\lim_{z\to z_{dom}}\frac{G(z_{max})-G(z)}{\sqrt{z_{dom}-z}}=z^{max}\sqrt{\frac{2}{-Z''_{max,1}(\tilde{y}^m)}}\frac{F(\tilde{y}^m)+F'(\tilde{y}^m)}{z^{max}-\tilde{y}^m}
\eeqlb
\item[(iii)]If  $z_{dom}=z^{max}=z^*$, then $z_{dom}$ is a pole of $G(z)$, and
\beqlb\label{s4-d-17}
\lim_{z\to z_{dom}}\sqrt{z_{dom}-z}G(z)=\frac{(\sum_{i=1}^2\sigma_i^2k_i^2)F(Y_{max,0}(z^{max}))z^{max}}
 {\sqrt{(\sum_{i=1}^2\sigma_i^2k_i^2)\sigma_3^2(z^{max}-y^{min})}}.
\eeqlb
\end{itemize}
\end{lem}
{\it Proof:} We  first prove  case (i). From Lemma \ref{lem-d-1}, in such a situation, $z_{dom}$ is the zero of $H_2(Y_{max,0}(z),z)=Y_{max,0}(z)-z$.  From \eqref{s2-87}, we get that
\beqlb\label{s4-d-19}
Y_{max,1}(z^*)-z^*\neq 0.
\eeqlb
From \eqref{s4-d-18},
\beqlb\label{s4-d-20}
f(z)=\frac{2(z^*-z)z}{\sum_{i=1}^3\sigma_i^2k_i^2}.
\eeqlb
From \eqref{s4-d-9} and \eqref{s4-d-20}, we obtain that
\beqlb\label{s4-d-21}
G(z)&&=\frac{2z(Y_{max,1}(z)-z)}{f(z)\sum_{i=1}^3\sigma_i^2k_i^2}F(Y_{max,0}(z))=\frac{2(Y_{max,1}(z)-z)}{\sum_{i=1}^3\sigma_i^2k_i^2(z^*-z)}F(Y_{max,0}(z)).
\eeqlb
\eqref{s4-d-15} follows from \eqref{s4-d-21}, since $z_{dom}=z^*$.

Next, we prove  case (ii).  From Lemma \ref{s4-lem3}, \eqref{s4-48} and \eqref{s4-d-13}, we get that $F(z)$ is analytic at the point $\tilde{y}^m$.  Hence,
\beqlb\label{s4-d-32}
F(Y_{max,0}(z))&&=F(\tilde{y}^m)+F'(\tilde{y}^m)(Y_{max,0}(z)-\tilde{y}^m)+o(|Y_{max,0}(z)-\tilde{y}^m|).
\eeqlb
It follows from \eqref{s4-d-9} that
\beqlb\label{s4-d-33}
G(z)&&=\frac{z}{(Y_{max,0}(z)-\tilde{y}^m+\tilde{y}^m-z)}F(Y_{max,0}(z))=\frac{z((Y_{max,0}(z)-\tilde{y}^m)-(\tilde{y}^m-z))}{(Y_{max,0}(z)-\tilde{y}^m)^2
-(\tilde{y}^m-z)^2}F(Y_{max,0}(z)).\qquad\quad
\eeqlb
From \eqref{s4-d-32} and \eqref{s4-d-33},
\beqlb\label{s4-45}
G(z)&&=\frac{F(\tilde{y}^m)z(Y_{max,0}(z)-\tilde{y}^m)}{(Y_{max,0}(z)-\tilde{y}^m)^2
-(\tilde{y}^m-z)^2}\nonumber
\\&&\qquad-(\tilde{y}^m-z)z\frac{F(\tilde{y}^m)+F'(\tilde{y}^m)(Y_{max,0}(z)-\tilde{y}^m)}{(Y_{max,0}(z)-\tilde{y}^m)^2-(\tilde{y}^m-z)^2}
+o(|Y_{max,0}(z)-\tilde{y}^m)|).
\eeqlb
Next, we consider the term $Y_{max,0}(z)-\tilde{y}^m$.  From Lemma \ref{lem-10a}  and \eqref{s4-d-13} that $Z_{max,1}(y)$ is analytic at the point $\tilde{y}^m$. Hence
\beqlb\label{s4-d-34}
Z_{max,1}(y)=Z_{max,1}(\tilde{y}^m)+Z'_{max,1}(\tilde{y}^m)(y-\tilde{y}^m)
+\frac{1}{2}Z''_{max,1}(\tilde{y}^m)(y-\tilde{y}^m)^2+o(|y-\tilde{y}^m|^2).\qquad
\eeqlb
Since $Z_{max,1}(y)$ takes the maximum at the point $\tilde{y}^m$ on $[y^{min},\;y^{max}]$,
\beqlb\label{s4-d-35}Z'_{max,1}(\tilde{y}^m)=0.\eeqlb
On the other hand, from \eqref{s5-2a} and  \eqref{s4-d-13},
\beqlb\label{s4-37}
Z_{max,1}(y)-Z_{max,1}(y^{m})<0.
\eeqlb
From \eqref{s4-d-34},  \eqref{s4-d-35} and \eqref{s4-37},
\beqlb\label{s4-d-36}
\tilde{y}^m-y=\sqrt{\frac{2(Z_{max,1}(y^{m})-Z_{max,1}(y))}{-Z''_{max,1}(\tilde{y}^m)}}+o(|y-\tilde{y}^m|).
\eeqlb
Similar to \eqref{s4-38},  we have that for $z$  close to $z^{max}$
\beqlb\label{s4-39}
Y_{max,0}(z)=y\;\textrm{and}\;Z_{max,1}(y)=z.
\eeqlb
Combing \eqref{s4-d-36} and \eqref{s4-39}, we obtain that
\beqnn
\tilde{y}^m-Y_{max,0}(z)=\sqrt{\frac{2(z^{max}-z)}{-Z''_{max,1}(\tilde{y}^m)}}+o(|z-z^{max}|^{\frac{1}{2}}).
\eeqnn
Hence,
\beqlb\label{s4-d-37}
\tilde{y}^m-Y_{max,0}(z)=(z^{max}-z)^{\frac{1}{2}}\sqrt{\frac{2}{-Z''_{max,1}(\tilde{y}^m)}}+o(|z-z^{max}|^{\frac{1}{2}}).
\eeqlb
From \eqref{s4-45} and \eqref{s4-d-37}, we obtain that
\beqlb\label{s4-46}
G(z)&&=(z^{max}-z)^{\frac{1}{2}}\frac{F(y^{m})z^{max}}{
-(\tilde{y}^m-z^{max})^2}\sqrt{\frac{2}{-Z''_{3,1}(\tilde{y}^m)}}+\frac{z^{max}F(\tilde{y}^m)}{(\tilde{y}^m-z^{max})}
\nonumber
\\&&+(z^{max}-z)^{\frac{1}{2}}\frac{F'(\tilde{y}^m)z^{max}\sqrt{\frac{2}{-Z''_{max,1}(\tilde{y}^m)}}}
{-(\tilde{y}^m-z^{max})}
+o(|z-\tilde{y}^m)|^{\frac{1}{2}}).
\eeqlb
Combining \eqref{s4-45} and \eqref{s4-d-37}, we obtain that
\beqlb\label{s4-d-38}
\lim_{z\to z^{max}}\frac{G(z)-G(z^{max})}{\sqrt{y^{max}-z}}=z^{max}\sqrt{\frac{2}{-Z''_{3,1}(\tilde{y}^m)}}\frac{F(\tilde{y}^m)+F'(\tilde{y}^m)}{z^{max}-\tilde{y}^m}.
\eeqlb

Finally, we prove  case (iii).
Due to Lemma \ref{lem-10}, we obtain that
\beqlb\label{s4-d-30}
H_2(Y_{max,0}(z^{max}),z^{max})=0.
\eeqlb
Hence,
\beqlb\label{s4-d-29}
&&\sqrt{z^{max}-z}\frac{z}{Y_{max,0}(z)-z}F(Y_{max,0}(z))\nonumber
\\&&\qquad\qquad=\frac{z\sqrt{z^{max}-z}}{H_2(Y_{max,0}(z),z)-H_2(Y_{max,0}(z^{max}),z^{max})}F(Y_{max,0}(z)).
\eeqlb
 From \eqref{s4-48}, \eqref{s4-d-13} and \eqref{s4-d-30}, $F(z)$ is analytic at $\tilde{y}^m=Y_{max,0}(y^{max})$. Therefore,
 \beqlb\label{s4-d-31}
 \lim_{z\to z^{max}}\sqrt{z^{max}-z}\frac {z}{Y_{max,0}(z)-z}F(Y_{max,0}(z))
 =\frac{(\sum_{i=1}^2\sigma_i^2k_i^2)F(Y_{max,0}(z^{max}))z^{max}}
 {\sqrt{(\sum_{i=1}^2\sigma_i^2k_i^2)\sigma_3^2(z^{max}-y^{min})}}.
 \eeqlb
\qed

We are now in the position to obtain  asymptotic properties of $\phi_2(0,0,z)$ and $\hat{\phi}_2(0,0,z)$ at the dominant singularity $z_{dom}$.
\begin{lem}\label{lem-d-4}
For the asymptotic behaviors of $\phi_2(0,0,z)$ and $\hat{\phi}_2(0,0,z)$ around the dominant singularity $z_{dom}$, we have
\begin{description}
\item[(i)] If $z_{dom}=z^*<z^{max}$, then
\beqlb\label{s4-d-22}
\lim_{z\to z^*}(z^*-z)\phi_2(0,0,z)&&=C_1(z^*),
\\ \label{s4-d-23}
\lim_{z\to z^*}(z^*-z)\hat{\phi}_2(0,0,z)&&=C_2(z^*);
\eeqlb
\item[(ii)] If $z_{dom}=z^*=z^{max}$, then
\beqlb\label{s4-d-22-1}
\lim_{z\to z^*}\frac{\phi_2(0,0,z^{*})-\phi_2(0,0,z)}{\sqrt{z^*-z}}&&=C_3(z^*),
\\ \label{s4-d-23-1}
\lim_{z\to z^*}\frac{\hat{\phi}_2(0,0,z^*)-\hat{\phi}_2(0,0,z)}{\sqrt{(z^*-z)}}&&=C_4(z^*);
\eeqlb
\item[(iii)] If $z_{dom}=z^{max}<z^*$, then
\beqlb\label{s4-d-22-2}
\lim_{z\to z^*}\sqrt{(z^*-z)}\phi_2(0,0,z)&&=C_5(z^*),
\\ \label{s4-d-23-3}
\lim_{z\to z^*}\sqrt{(z^*-z)}\hat{\phi}_2(0,0,z)&&=C_6(z^*).
\eeqlb
\end{description}
Here $C_i(z^*),i=1,\cdots,6$, are non-zero constants.
\end{lem}
{\it Proof:} Here, we only prove  case (i), other cases can be proved in the same fasion. It follows from \eqref{s4-97} that we only need focus on $z\in\big(2\frac{(c_3-\lambda_3-c_2)}{\sigma_3}, z^{max}\big)$. From \eqref{s2-87}, we get that
\beqlb\label{s6-8}
Y_{max}(z)\geq 0\;\textrm{for all}\;z\in \big(2\frac{(c_3-\lambda_3-c_2)}{\sigma_3}, z^{max}\big).
\eeqlb
Combing \eqref{s2-10} and \eqref{s6-8}, we get
\beqlb\label{s6-11}
\hat{\phi}_2(0,0,z)\geq \phi_2(0,0,z)
\eeqlb
for any $z\in\big(2\frac{(c_3-\lambda_3-c_2)}{\sigma_3}, z^{max}\big)$.
If  case (i) would not hold, then, from Lemmas \ref{s4-lem1}, \ref{lem-d-2} and \ref{lem-d-3}, we should have
\beqlb\label{s6-12}
C_1(z^*)=C_2(z^*)=\infty.
\eeqlb
If $k_1>1$, then from \eqref {s6-11} we have
\beqlb\label{s6-13}
G(z)=\phi_2(0,0,z)+\frac{1}{k_1-1}\hat{\phi}_2(0,0,z)\geq \frac{k_1}{k_1-1}\phi_2(0,0,z).
\eeqlb
From \eqref{s6-12} and \eqref{s6-13}, we get that
\beqlb\label{s6-15}
\lim_{z\to z^*}(z^*-z)G(z)=\infty,
\eeqlb
which contradicts to Lemma \ref{lem-d-3}.

On the other hand, if $0<k_1<1$, then from \eqref {s6-11}, we have
\beqlb\label{s6-16}
G(z)=\phi_2(0,0,z)+\frac{1}{k_1-1}\hat{\phi}_2(0,0,z)\leq \frac{k_1}{k_1-1}\hat{\phi}_2(0,0,z).
\eeqlb
Under this,  it is easily to check that
\beqlb\label{s6-17}
G(z)<0\;\textrm{for}\;z\in \big(2\frac{(c_3-\lambda_3-c_2)}{\sigma_3}, z^{max}\big).
\eeqlb
Hence, from Lemma \ref{lem-d-3}, we get that
\beqlb\label{s6-18}
-\infty<\lim_{z\to z^*}(z^*-z)G(z)<0.
\eeqlb
However, from \eqref{s6-11} and \eqref{s6-16}, we have
\beqlb\label{s6-19}
\lim_{z\to z^*}(z^*-z)G(z)=-\infty,
\eeqlb
which contradicts to \eqref{s6-18}. From above arguments, case (i) is proved.

Now we show that $C_i(z^*)$, $i=1,2$ are non-zero. It follows from \eqref{s4-d-15}, \eqref{s6-11} and \eqref{s6-13} that $C_2(z^*)\neq 0$. Now we assume that $C_1(z^*)=0$. Then from  \eqref{s2-10}, we have
\beqlb\label{4-2}
\hat{\phi}_2(0,0,z)&&=\int_{\R_{+}^{3}}\exp\big\{k_1Y_{max,0}(z)x_1+zx_3\big\}V_2(dx)\nonumber
\\&&<\frac{1}{2}\Big(\int_{\R_{+}^{3}}\exp\big\{2k_1Y_{max,0}(z)x_1\big\}V_2(dx)+\int_{\R_{+}^{3}}\exp\big\{2zx_3\big\}V_2(dx)\Big).
\eeqlb
Hence, as $z\to z^*$, from \eqref{4-2} we have
\beqlb\label{4-3}
C_2(z^*)<\tilde{C}_2(z^*)+C_1(z^*),
\eeqlb
where $$\tilde{C}_2(z^*)=\frac{1}{2}(z^*-z)\int_{\R_{+}^{3}}\exp\big\{2k_1Y_{max,0}(z)x_1\big\}V_2(dx).$$
On the other hand, it is obvious that
\beqlb\label{4-4}
\frac{1}{2}\int_{\R_{+}^{3}}\exp\big\{2k_1Y_{max,0}(z)x_1\big\}V_2(dx)<\frac{1}{2}\int_{\R_{+}^{3}}\exp\big\{2k_1Y_{max,0}(z)x_1+2zx_3\big\}V_2(dx).
\eeqlb
Finally, we note that
\beqlb\label{4-6}
Y_{max,0}(z)\to z,\; \textrm{as}\; z\to z^*.
\eeqlb
Letting $z\to z^*$ yields
\beqlb\label{4-5}
\tilde{C}_2(z^*)<C_2(z^*).
\eeqlb
Hence, \eqref{4-3} and \eqref{4-5} contradict to $C_1(z^*)=0$. From above arguments, we get that $C_i(z^*)>0$, $i=1,2$.

\qed

By Lemma \ref{s2-lem1},  in order to get tail asymptotics for the marginal $L_3$, we need to focus on $\phi_2(0,0,z)$. From asymptotic properties of  $\phi_2(0,0,z)$ obtained above, we can apply the Tauberian-like Theorem given below to transform asymptotic properties of  $\phi_2(0,0,z)$ to that of the marginal distribution $\P(L_3\leq x)$. To use the Tauberian-like theorem, we need study some properties of  $\phi_2(0,0,z)$ around the point $z_{dom}$. By Lemma \ref{s4-lem1}, there is exactly one dominant singularity for $\phi_2(0,0,z)$.  By Lemma \ref{lem-d-1}, there are two candidates for the dominant singularity $z_{dom}$ of $\phi_2(0,0,z)$:
\begin{description}
\item[(1)] A pole, i.e., a zero of $Y_{max,0}(z)-z$; or
\item[(2)] branch point $z^{max}$.
\end{description}

For each of  these two cases, we show that the analytic condition of the unknown function satisfies the Tauberian-like theorem.
\begin{lem}\label{lem-7}
If $z_{dom}<z^{max}$, then  there exists an $\epsilon>0$ such that $\phi_2(0,0,z)$ is analytic for $ \Re (z)< z_{dom}+\epsilon $ except for $z=z_{dom}$ and for each $a>0$
\beqlb\label{s2-62}
\sup_{\begin{array}{c}
z\notin B_{a}(z_{dom}) \\
\Re(z)< z_{dom}+\epsilon
\\
\end{array}} |\phi_2(0,0,z)|<\infty,\eeqlb
where $B_{a}(z_{dom})=\{z\in\mathbb{C}:|z-z_{dom}|<a\}$.
\end{lem}

\noindent{\it Proof:} From Lemma \ref{lem-d-1}, we see that if $z_{dom}<z^{max}$, then $z_{dom}$ is a pole of the function $\phi_2(0,0,z)$.  Hence,  $\phi_2(0,0,z)$ is analytic for $ \Re (z)< z_{dom}+\epsilon $ except for $z=z_{dom}$. It remains to show (\ref{s2-62}) for each $a>0$.  In such a case, $z_{dom}$ is a pole of $\phi_2(0,0,z)$.  It follows from Lemma \ref{lem-d-1} that $z_{dom}$ is a zero of $Y_{max,0}(z)-z$.  So
\beqlb\label{s2-64}
\sup_{\begin{array}{c}
       z\notin B_{a}(z_{dom}) \\ Re(z)< z_{dom}+\epsilon\\
      \end{array}
}\Big|\frac{1}{Y_{max,0}(z)-z}\Big|<\infty.
\eeqlb
From  \eqref{s4-d-14} and \eqref{s4-d-21}, we have for $\Re(z)<z_{dom}+\epsilon$
\beqlb\label{s2-65}
F\big(Y_{max,0}(z)\big)<\infty.
\eeqlb
Finally, we can easily get that
\beqlb\label{s2-66}
k_1-1<\infty.
\eeqlb
Equations (\ref{s2-64}) to (\ref{s2-66}) yield (\ref{s2-62}).  The proof is completed.
\qed
\begin{lem}\label{lem-8}
If $z_{dom}=z^{max}$, then $\phi(0,0,z)$ is analytic in $\mathbb{G}_{\delta_0}(z^{max})$. Moreover, for each $a>0$,
\beqnn
 \sup_{\begin{array}{c}
 z\in\mathbb{G}_{\delta_0}(z^{max})\\ z\notin B_{a}(z^{max}) \\
 \end{array}}|\phi_2(0,0,z)|<\infty.
\eeqnn
\end{lem}

\noindent{\it Proof:}  We first show that $\phi_2(0,0,z)$ is analytic on $z\in \mathbb{G}_{0}(z^{max})$. It follows from Lemma \ref{s4-lem1} that $\phi_2(0,0,z)$ is analytic for $\Re(z)<z_{dom}$. Furthermore, by \eqref{s4-69}, we have $z_{dom}>0$.  Hence, in order to prove the lemma, it suffices to show that
$\phi_2(0,0,z)$ is analytic on $ z\in \mathbb{G}_{\delta_0}(z^{max})\cap\{z\in\C: \Re(z)>0\}$.

Since $z_{dom}=z^{max}$, from Lemma \ref{lem-d-1}, we must have $z^*\geq z^{max}$. We first assume that
\beqlb\label{s6-10}z^{max}<z^*.
\eeqlb
Combing \eqref{s4-d-13} and Lemma \ref{s4-lem3}, we have that $F(Y_{max,0}(z))$ is analytic at $\mathbb{G}_{\delta_0}(z^{max})$ .   Hence, from \eqref{s4-d-9} and Corollary \ref{s3-cor1}, we can get the lemma.

Next, we assume that $z^{max}=z^*$.  The proof of this case is the combination of the proof of Lemma \ref{lem-7} and that of the  case \eqref{s6-10}.  So, we omit the details of the proof here.\qed

\section{Exact Tail Asymptotics for Marginal Distributions}

From the arguments in the previous section,  the asymptotic behavior and properties of $\phi_2(0,0,z)$ around the point $z_{dom}$ have been obtained.  In this section, we apply these results to get exact tail asymptotics for the marginal $L_3$.  Here, we also note that asymptotic behaviour of the marginal $L_3$ is closely related to the two points $z^*$ and $z^{max}$, which are  the candidates for $z_{dom}$. In practice,  we need to determine which one should be chosen as the dominant singularity $z_{dom}$.  In fact, we have the following lemma.
\begin{lem}\label{lem-1a}
$z^*$  exists between $(0,\;z^{max}]$ if and only if $Y_{max,0}(z^{max})\geq z^{max}$.
\end{lem}
{\it Proof:}  If $z^{max}=z^*$, one can easily see that the lemma holds. Next, we assume $z^{max}\neq z^*$. From \eqref{s2-87}, we obtain that $Y_{max,0}(z)$ is increasing on $(\frac{2(c_3-\lambda_3-c_2)}{\sigma_3^2},\;z^{max}].$   We first assume that $z^*$ exists in  $(0,\;z^{max})$. Since
\beqnn
0<z^*=Y_{max,0}(z^{*}),
\eeqnn
we  have
\beqlb\label{s5-7}
z^*>\frac{2(c_3-\lambda_3-c_2)}{\sigma_3^2}.
\eeqlb
Therefore
\beqlb\label{s5-4}
\tilde{y}^{max}:=Y_{max,0}(z^{max})>Y_{max,0}(z^{*}).
\eeqlb
On the other hand, we note that the line $H_2(y,z)=z-y=0$ intersects the ellipse $H(ky,y,z)=0$ at one point except for the point $(0,\;0)'$.  From \eqref{s2-97}, we know that the point $(Y_{max,0}(z^*),z^*)'$  is the other intersection point of $H_2(y,z)=0$ and $H(ky,y,z)=0$.  Hence, we must have
\beqlb\label{s5-5}
\tilde{y}^{max}>z^{max}.
\eeqlb

Next, we assume
\beqlb\label{s5-6}
\tilde{y}^{max}> z^{max}.
\eeqlb
We prove that $z^*$ belongs to $(0,\;z^{max})$.  From \eqref{s5-6}, we obtain that the point $(\tilde{y}^{max},z^{max})'$ is above the line $H_2(y,z)=0$.  From \eqref{s2-87}, we get that the point $\Big(Y_{max,0}(\frac{2(c_3-\lambda_3-c_2)}{\sigma_3^2}),\frac{2(c_3-\lambda_3-c_2)}{\sigma_3^2}\Big)'$ is below the line $H_2(y,z)=0$.  On the other hand,  $Y_{max,0}(z)$ is continuous on $\big(\frac{2(c_3-\lambda_3-c_2)}{\sigma_3^2},\;z^{max}\big)$. By the above arguments, one can get that the lemma holds.
\qed
\begin{rem}
From the proof of Lemma \ref{lem-1a}, we can get that if  $z^*$ exists, it is unique.
\end{rem}

 Once the dominant singularity is determined, we need to evaluate it.  In fact,
from Lemma \ref{lem-a1}, we can get $z^*$. On the other hand,  from Lemma \ref{lem-10}, we can obtain the value of $z^{max}$. Actually, we have

\beqlb\label{s6-14}
z^{max}=&&\frac{c_3-\lambda_3-c_2}{\sigma_3^2}+\nonumber
\\&&\frac{\sqrt{\big(\sigma_1^2+\sigma_2^2\big)^2
(c_3-\lambda_3-c_2)^2+\sigma_3^2(\sigma_1k^2
+\sigma_2^2)\big((c_1-\lambda_1)k_1
+c_2-\lambda_2-c_1\big)^2}}{\sigma_3^2(\sigma_1^2k^2+\sigma_2^2)}.
\eeqlb

After finding the values of $z^*$ and $z^{max}$, we can discuss the tail asymptotics of the marginal distribution. From Lemma \ref{s2-lem1},  we can see that the asymptotic behavior of $\phi(0,0,z)$ is closely related to that of $\phi_2(0,0,z)$. So, we first state the asymptotic behaviors of $\phi_2(0,0,z)$. The following lemma follows from Lemmas \ref{lem-d-3} and \ref{lem-d-4}.
\begin{lem}\label{lem-12}
For the function $\phi_2(0,0,z)$, a total of three types asymptotics exist as $x$ approaches to $z_{dom}$, based on the detailed properties of $z_{dom}$.
\begin{itemize}
\item[Case 1:]  If $z_{dom}=z^*<z^{max}$, then
    \beqlb\label{s6-1a}
    \lim_{z\to z_{dom}}(z_{dom}-z)\phi_2(0,0,x)=K_1(z_{dom});
   \eeqlb
\item[Case 2:] If $z_{dom}=z^*=z^{max}$, then
\beqlb\label{s6-2a}
 \lim_{z\to z_{dom}}\sqrt{z_{dom}-z}\phi_2(0,0,x)=K_2(z_{dom});
\eeqlb
\item[Case 3:] If $z_{dom}=z^{max}<z^*$, then
\beqlb\label{s6-3a}
 \lim_{z\to z_{dom}}\frac{\phi_2(0,0,x)-\phi_2(0,0,z_{dom})}{\sqrt{z_{dom}-z}}=K_3(z_{dom});
\eeqlb
\end{itemize}
where $K_i(z_{dom})$, $i=1,2,3$, are non-zero constants depending on the point $z_{dom}$.
\end{lem}
From Lemmas \ref{s2-lem1} and \ref{lem-12}, we can easily obtain asymptotic behavior of $\phi(0,0,z)$. In fact, we have:
\begin{lem}\label{lem-13}
For the moment generating function $\phi(0,0,x)$, a total of three types asymptotics exist as $z$ approaches to $z_{dom}$, based on the detailed property of $z_{dom}$.
\begin{itemize}
\item[Case 1:]  If $z_{dom}=z^*<z^{max}$, then
    \beqlb\label{s6-1}
    \lim_{z\to z_{dom}}(z_{dom}-z)\phi(0,0,z)=\bar{K}_1(z_{dom})
   \eeqlb
\item[Case 2:] If $z_{dom}=z^*=z^{max}$, then
\beqlb\label{s6-2}
 \lim_{z\to z_{dom}}\sqrt{z_{dom}-z}\phi(0,0,z)=\bar{K}_2(z_{dom})
\eeqlb
\item[Case 3:] If $z_{dom}=z_{max}<z^*$, then
\beqlb\label{s6-3}
 \lim_{z\to z_{dom}}\frac{\phi(0,0,z)-\phi(0,0,z_{dom})}{\sqrt{z^{dom}-z}}=\bar{K}_3(z_{dom})
\eeqlb
\end{itemize}
where $\bar{K}_i(z_{dom})$, $i=1,2,3,$ are non-zero constants depending on $z_{dom}$.
\end{lem}

Before we present the main result of this section, we need  the following technical tool.
Let $g(s)$ be the Laplace-transformation of $f(s)$, i.e,
\beqnn g(s)=\int_{0}^\infty e^{st}f(t)dt.
\eeqnn
Then, $g(s)$ is
analytic on the left half-plane. The singularities of $g(s)$  are
all in the right half-plane.  We have the following Tauberian-like theorem, which is due to Dai, Dawson and Zhao \cite{DDZ}.
\begin{thm}\label{TLT}
Assume that $g(z)$ satisfies the following conditions:
\begin{itemize}
\item[(1)]The left-most singularity of $g(z)$ is $\alpha_0$ with
$\alpha_0>0$. Furthermore, we assume that as $z\to \alpha_0$,
$$g(z)\sim(\alpha_0-z)^{-\lambda}$$ for some
$\lambda\in\C\setminus\mathbb{Z}_{\leq 0}$; \item[(2)]$g(z)$ is
analytic on $\mathbb{G}_{\epsilon_0}(\alpha_0)$ for some
$\epsilon_0\in(0,\frac{\pi}{2}]$; \item[(3)] $g(z)$ is bounded on
$\mathbb{G}_{\epsilon_1}(\alpha_0)$ for some $\epsilon_1>0$.
\end{itemize}
Then, as $t\to\infty$,
\beqlb\label{R1-16}
f(t)\sim e^{-\alpha_0 t}\frac{t^{\lambda-1}}{\Gamma(\lambda)},
\eeqlb
where $\Gamma(\cdot)$ is the Gamma function.
\end{thm}
Now, we state the main result of this paper.
\begin{thm}\label{C-1}
For the  tail of the marginal distribution  $\P(L_3>z)=V_3(z,\;\infty)$, we have the following tail asymptotic properties for large $z$:
\begin{description}
\item[Case 1:] If $z_{dom}=z^*<z^{max}$,
then
    $$V_3\big(z,\infty\big)\sim C_1e^{-z^{dom} z};$$
\item[Case 2:] If $z_{dom}=z^*=z^{max}$, then $$V_3\big(z,\infty\big)\sim C_2e^{-z_{dom}
z}z^{-\frac{1}{2}};$$

\item[Case 3:] If $z_{dom}=z^{max}<z^*$, then
$$V_3\big(z,\infty\big)\sim C_3e^{-z_{dom}
z}z^{-\frac{3}{2}};$$
\end{description}
where $C_i$, $i=1,2,3,$ are non-zero constants.
\end{thm}

\noindent{\it Proof:}
Cases (1) and (2) are  direct consequences of  Lemmas  \ref{s2-lem1},  \ref{lem-7}, \ref{lem-13}, and Theorem \ref{TLT}.

Next, we prove case (3).  From \eqref{s6-3}, we have
\beqlb\label{s6-4}
 \lim_{z\to z_{dom}}\sqrt{z_{dom}-z}\frac{\phi(0,0,z)-\phi(0,0,z_{dom})}{z_{dom}-z}=\bar{K}_3(z_{dom}).
\eeqlb
 From Dai and Miyazawza \cite{DM2011}, we get that  $$\frac{\phi(0,0,z_{dom})-\phi(0,0,z)}{z_{dom}-z}$$ is the moment generating function of the density function \beqlb\label{s6-6}\bar{f}(x)=e^{-z_{dom}x}\int_x^\infty e^{z_{dom}u}f(u)du,\eeqlb
 where $f(z)$  is the density function of the marginal distribution $\P(L_3<x)$.  Therefore, from Theorem \ref{TLT} and \eqref{s6-4}, we have
 \beqlb\label{s6-5}
 \bar{f}(z)\sim \hat{K}(z_{dom}) z^{-\frac{1}{2}}e^{-z_{dom}z},
 \eeqlb
 where $\hat{K}$ is a constant depending on $z_{dom}$.

 From \eqref{s6-6} and \eqref{s6-5}, we obtain that
 \beqlb\label{s6-9}
 \int_x^\infty e^{z_{dom}u}f(u)du\sim \hat{K}(z_{dom}) x^{-\frac{1}{2}}.
 \eeqlb
 Taking derivatives at the both sides of \eqref{s6-9}, we obtain that
 \beqlb\label{s6-7}
 f(x)\sim K_1(z_{dom})e^{-z_{dom}x}x^{-\frac{3}{2}},
 \eeqlb
 where $K_1$ is a constant. From \eqref{s6-7}, we conclude that  case (3) holds.
 \qed

\section{Tail Behaviours of Joint Stationary Distributions}

In this section, we  study  tail behaviours of the joint stationary distribution $\pi$ by using extreme value theory.  Before we state our main result of this section, we first introduce the domain of attraction of some extreme value distribution function $G(\cdot)$.
\defn{ ({\it Domain of Attraction}) Assume that $\big\{X_n=(X_n^{(1)},\cdots,X_n^{(d)})'\big\}$ are i.i.d.  multivariate random vectors with common distribution $\tilde{F}(\cdot)$  and the marginal distributions $\tilde{F}_i(\cdot)$, $i=1,\cdots,d$. If there exist normalizing constants $a_n^{(i)}>0$ and $b_n^{(i)}\in\R$, $1\leq i\leq d$, $n\geq 1$ such that as $n\to\infty$
\beqnn
\P\Big\{\frac{M_n^{(i)}-b_n^{(i)}}{a_n^{(i)}}\leq x^{(i)},1\leq i\leq d\Big\}&&=\tilde{F}^n\Big(a_n^{(1)}x^{(1)}+b_n^{(1)},\cdots,a_n^{(d)}x^{(d)}+b_n^{(d)} \Big)
\\&&\to G(x^{(1)},\cdots,x^{(d)}),
\eeqnn
where the maximum $M_n^{(i)}=\bigvee_{k=1}^{n}X_k^{(i)}$ is the componentwise maxima, then we call the distribution function $G(\cdot)$ a multivariate extreme value distribution function, and
$\tilde{F}$ is in the domain of attraction of $G(\cdot)$. We denote this by $\tilde{F}\in D(G)$.

For convenience, we let $F(x,y,z)$ denote the joint stationary distribution of $\{L(t)\}$ and  $F_i$,$i=1,2,3$, denote the stationary distribution of the $i$-th buffer content process. Miyazawa and Rolski \cite{MR2009}obtained exact tail asymptotics for marginal distributions $F_i$, $i=1,2$.   From Dai and Miyazawa \cite{DM2011}, Dai, Dawson and Zhao \cite{DDZ} and Theorem 6.2,  we can easily get the following lemma.
\begin{lem}\label{7-lem1}
For any $i\in\{1,2,3\}$, we have
\beqlb\label{7-a1}
1-F_i(x)\sim C_i \exp\{-\alpha_i x\}x^{\mu_i},
\eeqlb
where $\alpha_i$ is the dominant singularity of the moment generating function of the marginal distribution  $F_i$, and $\mu_i\in\{0,-\frac{1}{2},-\frac{3}{2}\}$ is the corresponding decay rate.
\end{lem}

From Lemma \ref{7-lem1}, we can get that
\begin{lem}\label{7-lem2}For any $i\in\{1,2,3\}$, we have
\beqnn F_i(x)\in D\big(G_1(x)\big),\eeqnn
where \beqnn G_1(x)=\exp\{-e^{-x}\}.\eeqnn
\end{lem}

\noindent{\it Proof:}  It follows from \eqref{7-a1} that as $x\to\infty$
\beqlb\label{7-a2}
F'_i(x)\sim \alpha_i \exp\{-\alpha_i x\} x^{\mu_i}\;\textrm{and}\;F''_i(x)\sim -\alpha^2_i \exp\{-\alpha_i x\} x^{\mu_i}.
\eeqlb
It follows from the asymptotic equivalence \eqref{7-a2} that
\beqlb\label{7-a3}
\lim_{x\to\infty}\frac{F''_1(x)\big(1-F_i(x)\big)}{\Big(F'_i(x)\Big)^2}=-1.
\eeqlb
Then, it follows from Proposition 1.1 in Resinck \cite[P.40]{R1987} that $F_i\in D(G_1)$. \qed

In the previous section, we obtained exact tail asymptotic properties of the marginal distributions. Now, based on these results, we can study the tail dependencies of joint stationary distributions. Before we state tail dependent result for  the stationary distributions $F(\cdot)$, we  introduce a technical lemma.

\begin{lem}\label{7-lem3}
Suppose that $\big\{X_n=(X_n^{(1)},X_n^{(2)},X_n^{(3)})'\big\}_{n\in\mathbb{N}}$ are i.i.d. random vectors in $\R^3$ with a common joint continuous distribution $\tilde{F}(\cdot)$, and the marginal distributions $\tilde{F}_i(\cdot)$, $i=1,2,3$. Moreover, we assume that $\tilde{F}_i(\cdot)$, $i=1,2,3$, are all in the domain of attraction of some univariate extreme value distribution $\hat{G}_1(\cdot)$, i.e., there exist $a_n^{(i)}$ and $b_n^{(i)}$ such that as $n\to\infty$
\beqnn
\tilde{F}^n_i\Big(a_n^{(i)}x+b_n^{(i)}\Big)\to \hat{G}_1(x),
\eeqnn
then, the following are equivalent.
\begin{itemize}
\item[(1)] $\tilde{F}$ is in the domain of attraction of a product measure, that is,
\beqnn
\tilde{F}^n\Big(a_n^{(i)}x^{(i)}+b_n^{(i)},i=1,2,3\Big)\to \Pi_{i=1}^3 \hat{G}_1\big(x^{(i)}\big).
\eeqnn
\item[(2)] For any  $1\leq i<j\leq 3$,
\beqlb\label{7-a39}
\lim_{t\to\infty}\P\Big(X^{(i)}>t, X^{(j)}>t\Big)/\big(1-\tilde{F}_q(t)\big)\to 0,
\eeqlb
where $q\in\{i,j\}$.
\end{itemize}
\end{lem}
By a slight modification of the proof of Proposition 5.27  in Rensick \cite[P.296]{R1987}, we can prove the above lemma.  Hence, we omit the detail here.

For the joint stationary distribution $F$, we have the following tail dependence.
\begin{lem}\label{7-lem4}
The joint stationary distribution function $F(\cdot)$ is asymptotically  independent, that is, there exist
 $a_n(\mu_i,\alpha_i)$ and $b_n(\mu_i,\alpha_i)$, $i=1,2,3$, such that
\beqnn
F^n(a_n(\mu_i,\alpha_i) x^{(i)}+b_n(\mu_i,\alpha_i),i=1,2,3 )\to \Pi_{i=1}^3 G_1(x^{(i)}),\;\text{as}\; n\to\infty.
\eeqnn
\end{lem}

\noindent{\it Proof:}
 Withou loss of generality, we assume that $L(0)=0$. We prove \eqref{7-a39}, an equivalent statement.  Here, we let $i=1$ and $j=2$ for simplicity. Other cases can be proved in the same fashion. We construct a new process such that the stationary tail probability is an  upper bound of the stationary tail probability $\P\{L\geq z\}$. Let $\hat{X}(t)=B(t)+\Lambda t$, which is  a three-dimensional Brownian motion, and
\beqnn
\hat{Y}(t)=-[R^{-1}\hat{X}(t)\wedge R^{-1}\Lambda t].
\eeqnn
Then, from Konstantopoulos, Last and Lin \cite[Proposition 1]{KLL2004}, we get that  for any $z=(z_1,z_2,z_3)'\in\R^3_+$
\beqlb\label{7-a4}
\P\{L(t)\geq z\}\leq \P\{\hat{L}(t)\geq z\},
\eeqlb
where the operations are performed component-wise, and
\beqnn
\hat{L}(t)= \hat{X}(t)+R\hat{Y}(t).
\eeqnn
On the other hand, for any $z\in\R_+^3$
\beqlb\label{7-a32}
\P\{L\geq z\}=\lim_{t\to\infty}\P\big\{L(t)\geq z\big\}=\liminf_{t\to\infty}\P\big\{L(t)\geq z\big\}\leq \P\big\{L(1)\geq z\big\}
\eeqlb
For any $\tilde{z}=\{\tilde{z}_1,\tilde{z}_2)'\in\R_+^2$, let
\beqnn
\bar{F}_{12}(\tilde{z})=\P\big\{L_1\geq \tilde{z}_1, L_2\geq \tilde{z}_2\}.
\eeqnn
Then, from \eqref{7-a32},
\beqlb\label{7-a5}
\bar{F}_{12}(\tilde{z})\leq\P\{\tilde{L}(1)\geq\tilde{z}\},
\eeqlb
where $\tilde{L}(t)=AL(t)$ with $$A=\begin{bmatrix}1\;&0\;&0\\0&1&0\end{bmatrix}.$$

From \eqref{7-a4} and \eqref{7-a5}, we get that
\beqlb\label{7-a6}
\bar{F}_{12}(\tilde{z})\leq  \P\{\hat{L}(1)\geq \tilde{z}\}
\eeqlb
Hence, for any $\tilde{z}=(\tilde{z}_1,\tilde{z}_2)'\in\R^3_+$,
\beqlb\label{7-a7}
\P\{\tilde{L}(1)\geq \tilde{z}\}&&\leq \P\{A\big(\hat{X}(1)-\Lambda\big)\geq\tilde{z}\}.
\eeqlb

It is obvious  that $A\big(\hat{X}(1)-\Lambda\big)$ is bivariate Gaussian vector with the correlation coefficient being less than one. On the other hand, from \eqref{7-a7}, we have that for large enough $z\in\R_+$
\beqlb\label{7-a9}
\limsup_{z\to\infty}\frac{\bar{F}_{12}(z,z)}{\bar{F}_1(z)}\leq \limsup_{z\to\infty}\frac{ \P\{\hat{L}(1)\geq (z,z)'\}}{\bar{F}_1(z)}.
\eeqlb
Finally, we get that
\beqlb\label{4-a10}
&&\limsup_{z\to\infty}\frac{\P\{\hat{L}_1(1)\geq z,\hat{L}_2(1)\geq z\}}{\P\{L_1 \geq z\}}\nonumber
\\
&&\hspace{2cm}=\limsup_{z\to\infty}\frac{\P\{\hat{L}_1(1)\geq z, \hat{L}_2(1)\geq z\}}{\P\{\hat{L}_1(1)\geq z\}}\frac{\P\{\hat{L}_1(1)\geq z\}}{\P\{L_1\geq z\}}\nonumber
\\&&\hspace{2cm}\leq \limsup_{z\to\infty} \frac{\P\{\hat{L}_1(1)\geq z, \hat{L}_2(1)\geq z\}}{\P\{\hat{L}_1 \geq z\}}=0,
\eeqlb
where the  inequality is obtained by using
\beqnn
\P\{\hat{L}_1(1)\geq z\}/\P\{L_1\geq z\}\to 0,\;\text{as}\; z\to\infty.
\eeqnn
From above arguments, we get that
\beqlb\label{7-a11}
\lim_{z\to\infty}\frac{\bar{F}_{12}(z,z)}{\bar{F}_1(z)}=0.
\eeqlb
From \eqref{7-a11} and Lemma \ref{7-lem3}, we get the lemma.
\qed
\begin{rem}
For $a_n(\mu_i,\alpha_i)$ and $b_n(\mu_i,\alpha_i)$, $i=1,2,3$ in Lemma \ref{7-lem4}, we can use tail equivalence to get their explicit expressions.  Since they are not the focus of this paper, we will not elaborate it here.
\end{rem}
Now, we present the main result of this section.
\begin{thm}\label{thm} As $(x,y,z)'\to(\infty,\infty,\infty)'$,
\beqlb\label{7-a12}\P\big\{L_1\geq x, L_2\geq y, L_3\geq z\big\}/\Big(K x^{\mu_1}y^{\mu_2}z^{\mu_3}\exp\big\{-(\alpha_1x+\alpha_2y+\alpha_3z)\big\}\Big)\to 1,
\eeqlb
where $\alpha_i$ is the dominant singularity of $L_i$, and $\mu_i\in\{0,-\frac{1}{2},-\frac{3}{2}\}$ is the exponent corresponding to $\alpha_i$ in Lemma \ref{7-lem1}
\end{thm}

\noindent{\it Proof:}
 To prove this theorem, we first need a transformation.
Let $\bar{X}=(\bar{X}_1,\bar{X}_2,\bar{X}_3)'$ be a random vector with the joint distribution $\tilde{F}(x,y,z)$  and marginal distributions $\tilde{F}_i(x)$, $i=1,2,3$. Then we make the following transformation:
\beqlb\label{7-a13}
X^*_i=\frac{-1}{\log \big(\tilde{F}_i(\bar{X}_i)\big)},\;\text{for}\;i=1,2,3.
\eeqlb
By the transformation  in \eqref{7-a13}, we transform each marginal $\bar{X}_i$ of a random vector $\bar{X}$  to a unit Fr\'echet variable $X_i^*$, that is,
\beqnn
\P\{X^*_i<x\}=\exp\{-\frac{1}{x}\}\;\text{for}\;x\in\R_+.
\eeqnn
Hence, for the trivariate extreme value distribution $G(x,y,z)$
\beqlb\label{7-a14}
G^*(x,y,z)=G\Bigg(\Big(\frac{-1}{\log G_1}\Big)^{-1}\big(x\big),\Big(\frac{-1}{\log G_1}\Big)^{-1}\big(y\big),\Big(\frac{-1}{\log G_1}\Big)^{-1}\big(z\big)\Bigg),
\eeqlb
where $G^*(\cdot,\cdot,\cdot)$ is the joint distribution function  with the common marginal Fr\'echnet distribution $\Phi(x)=\exp\{-x^{-1}\}$.
Furthermore, for the stationary random vector $L$, define
\beqlb\label{7-a15}
Y_i=\frac{1}{1-F_i(L_i)},
\eeqlb
where $F_i$ is the marginal distribution of $L_i$.
Let $F^*(y_1,y_2,y_3)$ be the joint distribution function of $Y=(Y_1,Y_2,Y_3)'$. Then, it follows from Proposition 5.10 in Resnick \cite{R1987} and Lemma \ref{7-lem4} that
\beqlb\label{7-a16}
F^*(y_1,y_2,y_3)\in D\big(G^*(y_1,y_2,y_3)\big).
\eeqlb
By \eqref{7-a16}, we have that for any $y=(y_1,y_2,y_3)'\in\R_+^3$, as $n\to\infty$
\beqlb\label{7-a18}
(F^*(ny))^n\to G^*(y) .
\eeqlb
It follows from \eqref{7-a18} that
\beqnn
F^*(ny)\sim \big(G^*(y)\big)^{\frac{1}{n}}.
\eeqnn
By a simple monotonicity argument, we can replace $n$ in the above equation by $t$.
Then we have that as $t\to\infty$,
\beqlb\label{7-a19}
F^*(ty)\sim \big(G^*(y)\big)^{\frac{1}{t}}.
\eeqlb
On the other hand, by Lemma \ref{7-lem2}, for any $y\in\R_+$
\beqlb\label{7-a20}
F^*_i(ty)\sim \big(G^*_1(y)\big)^{\frac{1}{t}},\;\text{for any}\;i=1,2,3.
\eeqlb
Combing \eqref{7-a19} and \eqref{7-a20}, we get that as $t\to\infty$
\beqlb\label{7-a21}
F^*(ty)\sim F_1^*(ty_1)\cdot F_2^*(ty_2)\cdot F_3^*(ty_3).
\eeqlb
Let $C(u_1,u_2,u_3)$ be the copula of the random vector $(Y_1,Y_2,Y_3)'$, i.e.,
\beqlb\label{7-a27}
C\big(F_1^*(x),F_2^*(y),F_3^*(z)\big)=F^*(x,y,z).
\eeqlb
Furthermore, let $\hat{C}(u_1,u_2,u_3)$ be the corresponding  survival copula of $C$. Then, we have (see, for example, equation (2.46) in Schmitz \cite{SV2003} ):
\beqlb\label{7-a41}
\hat{C}(u_1,u_2,u_3)&&=\sum_{i=1}^3 u_i+\sum_{1\leq i<j\leq 3}C_{i,j}(1-u_i,1-u_i)\nonumber
\\&&\hspace{0.5cm}-C(1-u_1,1-u_2,1-u_3)-2.
\eeqlb
For convenience, for any $(x_1,x_2,x_3)'\in\R_+^3$, let $u_i(t)=\bar{F}^*_i(tx_i)$. Hence for any $t\in\R_+$,
\beqlb\label{7-a40}
\hat{C}\big(u_1(t),u_2(t),u_3(t)\big)=\bar{F}^*(tx_1,tx_2,tx_3),
\\ C\big(1-u_1(t),1-u_2(t),1-u_3(t)\big)=F^*(tx_1,tx_2,tx_3).\nonumber
\eeqlb
 Moreover, from \eqref{7-a21}, we get that as $t\to\infty$
\beqlb\label{7-a42}
C\big(1-u_1(t),1-u_2(t),1-u_3(t)\big)\sim \big(1-u_1(t)\big)\cdot\big(1-u_2(t)\big)\cdot\big(1-u_3(t)\big),
\eeqlb
and for any $1\leq i<j\leq 3$
\beqlb\label{7-a43}
C_{i,j}\big(1-u_i(t),1-u_j(t)\big)\sim \big(1-u_i(t)\big)\cdot \big(1-u_j(t)\big).
\eeqlb
From \eqref{7-a41}, \eqref{7-a42} and \eqref{7-a43}, we get that as $t\to\infty$
\beqlb\label{7-a24}
\hat{C}\big(u_1(t),u_2(t),u_3(t)\big)\sim  u_1(t)\cdot u_2(t)\cdot u_3(t),
\eeqlb
which is equivalent to for any $(x,y,z)'\in\R^3_+$
\beqlb\label{7-a25}
\lim_{t\to \infty}\frac{\bar{F}^*(tx,ty,tz)}{\bar{F}^*_1(tx)\cdot \bar{F}^*_2(ty)\cdot \bar{F}^*_3(tz)}=1.
\eeqlb
To prove our theorem, we need to show
\beqlb\label{7-a26}
\lim_{(x,y,z)'\to(\infty,\infty,\infty)'}\frac{\bar{F}^*(x,y,z)}{\bar{F}^*_1(x)\cdot \bar{F}^*_2(y)\cdot \bar{F}^*_3(z)}=1.
\eeqlb

From \eqref{7-a27}, to prove \eqref{7-a26}, we only need to show
\beqlb\label{7-a28}
\lim_{(u,v,w)'\to(0,0,0)'\;\textrm{and}\; (u_1,u_2,u_3)'\in I^3}\frac{\hat{C}(u_1,u_2,u_3)}{u_1u_2u_3}=1,
\eeqlb
where $I=[0,\;1]$.
Note that
\beqlb\label{7-a29}
\lim_{x\to 0}\frac{1-\exp\{-x\}}{x}=1.
\eeqlb
Hence from \eqref{7-a25}, we get that
\beqlb\label{7-7}
\lim_{t\to 0}\frac{\hat{C}(tu_1,tu_2,tu_3)}{t^3u_1u_2u_3}=1.
\eeqlb
Here we should point out that the limit \eqref{7-a28} has the form of $\frac{0}{0}$. Hence, we apply the multivariate L'h$\hat{o}$pital's rule (see Theorem 2.1 in  Lawlor \cite{L2012}) to prove it.  Without much effort, we can construct a multivariate differential function $\tilde{C}(u,v,w)$ such that
\beqnn
\hat{C}(u_1,u_2,u_3)=\tilde{C}(u_1,u_2,u_3)\;\textrm{for all}\; (u_1,u_2,u_3)'\in I^3,
\eeqnn
and
\beqnn
\tilde{C}(tu_1,tu_2,tu_3)\sim t^3 u_1u_2u_3,\;\textrm{as}\;t\to 0.
\eeqnn

Hence we only have to
\beqlb\label{7-8}
\lim_{(u_1,u_2,u_3)'\to(0,0,0)'\;\textrm{and}\; (u_1,u_2,u_3)'\in I^3}\frac{\hat{C}(u_1,u_2,u_3)}{wuv}=\lim_{(u_1,u_2,u_3)'\to(0,0,0)'}\frac{\tilde{C}(u_1,u_2,u_3)}{u_1u_2u_3}=1.
\eeqlb
Near the origin $(0,0.0)'$, the zero sets  of both $\tilde{C}(u_1,u_2,u_3)$ and $u_1u_2u_3$  consist of the hypersurfaces $u_1=0$, $u_2=0$ and $u_3=0$.  By the multivariate \'Lh$\hat{o}$pital's rule  (see Theorem 2.1 in  Lawlor \cite{L2012}), to prove \eqref{7-8}, we need to show that for each component $E_i$ of $\R^3\setminus\mathcal{C}$, where $\mathcal{C}=\{u_1=0\}\cup\{u_2=0\}\cup\{u_3=0\}$, we can find a vector $\vec{z}$, not tangent to $(0,0,0)'$ such that $D_{\vec{z}}(u_1u_2u_3)\neq 0$ on $E_i$ and
\beqnn
\lim_{(u_1,u_2,u_3)'\to (0,0,0)'\; \textrm{and}\; (u_1,u_2,u_3)'\in E_i} \frac{D_{\vec{z}}\tilde{C}(u_1,u_2,u_3)}{D_{\vec{z}}(u_1u_2u_3)}=1
\eeqnn
For the component $E_1$ bounded by the hypersurfaces of $\mathcal{H}_i=\{(u_1,u_2,u_3)':(u_1,u_2,u_3)'\in\R_+^3\;\textrm{and}\; u_i=0\}$,  $i=1,2,3$, choose, say $\vec{z}=(1,1,1)'$, then $z$ is not tangent to any hypersurfaces $u_i=0$, $i=1,2,3$  at the point $(0,0,0)'$.   Next, we take the limits along the direction $\vec{z}=(1,1,1)'$. It follows from \eqref{7-7} and \eqref{7-a29} that
\beqlb\label{7-a30}
\lim_{(u_1,u_2,u_3)'\to (0,0,0)'\;\textrm{and}\;(u_1,u_2,u_3)'\in E_1} \frac{D_{\vec{z}}\tilde{C}(u_1,u_2,u_3)}{D_{\vec{z}}(u_1u_2u_3)}=1.
\eeqlb
Similar to \eqref{7-a30}, for any other components $E_i$, $i=2,\cdots,8$, we can find a vector $\vec{z}$ such that $z$ is not tangent to any hypersurfaces $u_i=0$, $i=1,2,3$  at the point $(0,0,0)'$. Moreover, we have
\beqlb\label{7-a31}
\lim_{(u_1,u_2,u_3)'\to (0,0,0)'\;\textrm{and}\;(u_1,u_2,u_3)'\in E_i} \frac{D_{\vec{z}}\tilde{C}(u_1,u_2,u_3)}{D_{\vec{z}}(u_1u_2u_3)}=1.
\eeqlb

From \eqref{7-8} to \eqref{7-a31} and  Lawlor \cite{L2012} that
\beqlb\label{7-a33}
\lim_{(u,v,w)'\to (0,0,0)'} \frac{\hat{C}(u_1,u_2,u_3)}{u_1u_2u_3}=1.
\eeqlb

Finally, it follows from \eqref{7-a15} that for any $(x,y,z)'\in\R_+^3$,
\beqlb\label{7-a34}
\P\{L_1\geq x,L_2\geq y,L_3\geq z\}&&=\P\big\{F_1(L_1)\geq F_1(x),F_2(L_2)\geq F_2(y),F_3(L_3)\geq F_3(z)\big\}\nonumber
\\&&=\P\big\{Y_1\geq \frac{1}{1-F_1(x)},Y_2\geq \frac{1}{1-F_2(y)},Y_3\geq \frac{1}{1-F_3(z)}\big\}\nonumber
\\&&=F^*\Big(\frac{1}{\bar{F}_1(x)},\frac{1}{\bar{F}_2(y)},\frac{1}{\bar{F}_3(z)}\Big).
\eeqlb
Combining \eqref{7-a33} and \eqref{7-a34}, we get that as $ (x,y,z)'\to(\infty,\infty,\infty)'$
\beqlb\label{7-a35}
\P\{L_1\geq x,L_2\geq y,L_3\geq z\}/\Bigg( \bar{F}^*_1\bigg(\frac{1}{\bar{F}_1(x)}\bigg)\cdot  \bar{F}^*_2\bigg(\frac{1}{\bar{F}_2(y)}\bigg)\cdot \bar{F}^*_3\bigg(\frac{1}{\bar{F}_3(z)}\bigg)\Bigg)\to 1.
\eeqlb
By \eqref{7-a35} and \eqref{7-a29}, we get that as $(x,y,z)'\to (\infty,\infty,\infty)'$
\beqlb\label{7-a37}
\P\{L_1\geq x,L_2\geq y,L_3\geq Z\}/\Big( \bar{F}_1(x)\cdot \bar{F}_2(y)\cdot \bar{F}_3(z)\Big)\to 1.
\eeqlb
Finally, it follows from Lemma \ref{7-lem1}  and \eqref{7-a37} that
\beqlb\label{7-a38}
\bar{F}_i(x)\sim K x^{\mu_i}\exp\{-\alpha_i x\},\;i=1,2,3.
\eeqlb
From \eqref{7-a37} and \eqref{7-a38}, we get the theorem. \qed

\section{Concluding Remarks}

In the previous sections, we obtained tail asymptotic properties for $L_3$, see Theorem \ref{C-1}, and asymptotic independence for $L$, see Theorem \ref{thm}.
An immediate question is: Can we generalize our study to the model with a dimension higher than three? To answer this question, we first recall the key components in our analysis for the 3-dimensional model:
\begin{enumerate}
\item The fundamental form, or the functional equation satisfied by the (unknown) moment generating functions of the joint stationary distribution and boundary measures (the counterpart to the equation  in (\ref{s2-69})).

By using It'o formula, such a relationship can be obtained for the $n$-dimensional model.

\item The kernel method, including analytic continuation of the unknown moment generating functions and asymptotic analysis. 

This seems to be the main challenge. It is our conjecture that the counterpart analytic continuation property (to Lemma \ref{lem-a-4}) is still there for the $n$-dimensional case. If this is true, the asymptotic analysis on the dominant singularity should prevail.

\item Based on the asymptotic analysis of the unknown moment generating functions, the same Tauberian-like Theorem will lead to the tail asymptotic properties for the boundary measures and marginal distributions, the counterpart to Theorem \ref{C-1}.

\item Furthermore, by extreme value theory and copula, similar to Theorem \ref{thm}, asymptotic independence for joint stationary distributions can be obtained. 

\end{enumerate}

\vspace{5mm}

\noindent{\bf Acknowledgments} This work was supported in part by the National Natural Science Foundation of China (No.71671104), the Fostering Project of Dominant Discipline and Talent Team of Shandong Province Higher Education Institutions, and the National Science and Engineering Research Council (NSERC) of Canada.

\end{document}